# SINGULAR CONTROL WITH STATE CONSTRAINTS ON UNBOUNDED DOMAIN

By Rami Atar[1] and Amarjit Budhiraja[2]

*Technion–Israel Institute of Technology and University of North Carolina*

We study a class of stochastic control problems where a cost of the form

$$\mathbb{E} \int_{[0,\infty)} e^{-\beta s} [\ell(X_s)\, ds + h(Y_s^\circ)\, d|Y|_s] \tag{0.1}$$

is to be minimized over control processes $Y$ whose increments take values in a cone $\mathbb{Y}$ of $\mathbb{R}^p$, keeping the state process $X = x + B + GY$ in a cone $\mathbb{X}$ of $\mathbb{R}^k$, $k \le p$. Here, $x \in \mathbb{X}$, $B$ is a Brownian motion with drift $b$ and covariance $\Sigma$, $G$ is a fixed matrix, and $Y^\circ$ is the Radon–Nikodym derivative $dY/d|Y|$. Let $\mathcal{L} = -(1/2)\mathrm{trace}(\Sigma D^2) - b \cdot D$ where $D$ denotes the gradient. Solutions to the corresponding dynamic programming PDE,

$$[(\mathcal{L} + \beta)f - \ell] \vee \sup_{y \in \mathbb{Y}: |Gy|=1} [-Gy \cdot Df - h(y)] = 0, \tag{0.2}$$

on $\mathbb{X}^o$ are considered with a polynomial growth condition and are required to be supersolution up to the boundary (corresponding to a "state constraint" boundary condition on $\partial \mathbb{X}$). Under suitable conditions on the problem data, including continuity and nonnegativity of $\ell$ and $h$, and polynomial growth of $\ell$, our main result is the unique viscosity-sense solvability of the PDE by the control problem's value function in appropriate classes of functions. In some cases where uniqueness generally fails to hold in the class of functions that grow at most polynomially (e.g., when $h = 0$), our methods provide uniqueness within the class of functions that, in addition, have compact level sets. The results are new even in the following special cases: (1) The one-dimensional case $k = p = 1$, $\mathbb{X} = \mathbb{Y} = \mathbb{R}_+$; (2) The first-order case $\Sigma = 0$; (3) The case where $\ell$ and $h$ are linear. The proofs

Received May 2004; revised April 2005.
[1]Supported in part by the Israel Science Foundation Grant 126/02.
[2]Supported in part by the Army Research Office Grant W911NF-04-1-0230.

combine probabilistic arguments and viscosity solution methods. Our framework covers a wide range of diffusion control problems that arise from queueing networks in heavy traffic.

**1. Introduction.** This paper is concerned with a class of singular stochastic control problems with state constraints. The controlled diffusion process takes values in a closed convex cone $\mathbb{X}$. The cost is of the form (0.1), and the running cost, $\ell$, is not assumed to be bounded. The corresponding dynamic programming equation (0.2) is considered with a polynomial growth condition, and the role of a boundary condition on $\partial \mathbb{X}$ is played by the requirement that the solution be a supersolution up to the boundary. It is well known that classical solutions to *Hamilton–Jacobi–Bellman* (HJB)-type equations do not exist in general, and that an appropriate partial differential equations (PDE) framework is via the notion of viscosity solutions (cf. Section 2). Our main result is the characterization of the value function for the problem as the unique viscosity-sense solution of (0.2).

For an introduction to viscosity solutions and a list of relevant literature, the reader is referred to Crandall, Ishii and Lions [8]. Much of the motivation and first examples of this theory came from problems in optimal control. See [1] and [15] for accounts on viscosity solutions in deterministic and stochastic control and differential games. Control problems with state constraint were studied in the viscosity solution framework in the monograph by Lions [27] for the deterministic (first-order) case. The relation between control with state constraint and viscosity supersolutions of the corresponding PDE on the boundary was first observed and developed by Soner [32]. It was extended to more general first-order equations by Capuzzo-Dolcetta and Lions [6]. For stochastic control with state constraints (and no singular term), see [25, 26], the recent work of Ishii and Loreti [22] and references cited therein. See also [10] and references therein for models in mathematical finance with state constraints.

Many authors have contributed to the study of singular control of diffusions. One-dimensional problems were studied by direct analysis by Beneš, Shepp and Witsenhausen [2], Karatzas [23] and Harrison and Taksar [18]. HJB-type PDE associated with singular stochastic control were studied for their classical and weak (a.e.-sense) solutions by Evans [14], Menaldi and Robin [29], Chow, Menaldi and Robin [7] and Ishii and Koike [21]. Contributions to viscosity solutions for such HJB equations include [15, 28] and [31]. In [28] the authors consider a two-dimensional model that arises from the heavy traffic analysis of a queueing network and corresponds to a problem of the type studied in the current paper (with linear $\ell$ and $h = 0$). The authors use the viscosity solution framework to establish asymptotic optimality of proposed control schemes. The model studied in [31] arises from a problem in mathematical finance. The role of the singular control there is, in a sense, antipodal to its role in the current model: it is possible to force the state process to move to the origin and remain there, whereas in the current model it will be seen that the control can only contribute to a motion away from the origin [in the sense of equation (2.2)].

Some of the papers that address uniqueness of solutions to HJB equations on unbounded domains are as follows. The paper [29] mentioned above investigates a class of singular control problems (without state constraints and with $h = 0$) and gives various results on characterization of the value function as the maximal or the unique solution of the associated PDE in a suitable weak sense; however, viscosity solutions are not considered there. In [10] the authors study a drift control problem with state constraints and prove uniqueness within a class of concave functions. Results for a broad family of second-order degenerate elliptic PDE appear in [20] and [8]; however the results on unbounded domains therein do not cover PDEs associated with singular control or state constraints. Crandall and Lions [9] consider first-order equations on unbounded domain and, motivated by problems in optimal control, extend previous results on uniqueness within uniformly continuous functions (see references therein) to continuous sub-exponential functions. All of the results mentioned above fail to cover PDEs of the form (0.2). The delicate nature of the uniqueness problem in the current setting can be seen from several simple one-dimensional examples, that are provided in Section 2, where uniqueness fails. In particular, when $h = 0$, uniqueness may fail to hold in the class of functions that grow at most polynomially. In such cases, our approach establishes uniqueness within the class of functions that, in addition, have compact level sets.

A primary motivation for the problems considered in this paper comes from controlled queueing systems in heavy traffic. A formal diffusion approximation of such systems leads to a class of problems referred to as *Brownian control problems* (BCPs) (cf. [16, 17]). These, in turn can be transformed, using techniques introduced in [19], to singular control problems with state constraints, of the form considered in this paper. We will demonstrate that our result on the PDE characterization of the value function covers control problems arising from a broad family of stochastic networks (see Section 3).

We now remark on some of the key steps in the proof of our main result. There are two natural ways to define a value function for the problem: A strong formulation, in which infimum of the cost is taken over control processes $Y$ adapted to the Brownian motion, and a weak formulation, where the infimum ranges over all filtered probability spaces and all control processes adapted to the underlying filtration. Denoting the two resulting value functions by $V$ and $\overline{V}$, it is clear that $\overline{V} \leq V$ [see Section 2, equations (2.12) and (2.13) for precise definitions]. Solvability of the PDE by both $V$ and $\overline{V}$ is established by means of two different dynamic programming principles (DPPs) (Propositions 5.1 and 5.2). While the DPP associated with $V$ is quite standard and is essentially a consequence of the strong

Markov property of Brownian motion, the DPP for $\overline{V}$ is less straightforward. The latter relies on a representation of $\overline{V}$ as the infimum of the cost over controls that are, in an appropriate sense, of feedback form (cf. Section 5). Results of this type are well understood for absolutely continuous controls (cf. [4]), and go back to a fundamental martingale representation result due to Wong [33]. However, in presence of singular control and state constraint, our result on the DPP for $\overline{V}$ appears to be new. The proof also makes use of the so-called *Skorohod problem* to account for the state constraint. For a DPP for singular control problems (without state constraints), see [29].

The proof of uniqueness of solutions is carried out first for a mixed Dirichlet state-constraint boundary value problem on a bounded domain, and then lifted to the unbounded domain. Uniqueness on bounded domain uses tools from the theory of viscosity solutions, and although several key ingredients in the argument have been well developed in the literature, it appears that this paper is the first to prove uniqueness of viscosity solutions of the HJB equation for a stochastic singular control problem with state constraints. This treatment could, in fact, be carried out for a much more general second-order degenerate elliptic operator than $\mathcal{L}$. The limitation put on the operator comes from our treatment of the problem on unbounded domain, where certain estimates on the dynamics are used crucially [viz., (4.10)–(4.12)].

As already mentioned, the literature on second-order degenerate elliptic PDE on unbounded domains fails to capture uniqueness for equation (0.2). Our approach uses a verification argument that compares an arbitrary solution $u$ to the value function $V$. More precisely, it is first shown that $u$ solves a singular control problem on a bounded domain, with an exit cost equal to its value on the boundary, giving a variational representation for $u$ similar to a DPP (cf. Proposition 6.1). It is here that the uniqueness result on bounded domain is required. This sets the ground for comparing $V$ with $u$ by means of constructing an admissible control for one problem using the other and considering the control problem for $u$ on an increasing sequence of domains. In such a construction, a large time sub-exponential estimate on the controlled process is crucially used in obtaining the inequality $u \leq V$. For the inequality $V \leq u$, a control process for the problem associated with $V$ is constructed by suitably patching together a sequence of controls for the bounded domain problems.

The paper is organized as follows. Some notation is introduced at the end of this section. Section 2 introduces the control problem setting, the PDE, the main result and some examples of nonuniqueness. Section 3 demonstrates the applicability of the main result to problems that arise from queues in heavy traffic. Section 4 contains the bounded domain problem formulation and some preliminary lemmas. Section 5 proves solvability of the PDE by the value functions, based on the DPPs. Sections 6 and 7 establish uniqueness on bounded, and respectively, unbounded domain. In Section 8, the DPPs are proved. Finally, some auxiliary results are provided in the Appendix.

The following notation will be used.

For $\alpha \in \mathbb{R}^n$, $|\alpha|$ denotes the Euclidean norm.

For a set $S \subset \mathbb{R}^n$, $C^2(S)$ denotes the space of twice continuously differentiable functions on $S$. $C^c(S)$ denotes the class of continuous functions $f$ on $S$ for which all level sets $\{x \in S : f(x) \leq r\}$, $r \in \mathbb{R}$, are compact. $C_{\text{pol}}(S)$ denotes the class of continuous functions $f$ on $S$ for which there is a constant $a = a(f)$ such that $|f(x)| \leq a(1 + |x|)^a$, $x \in S$. $C_{\text{b}}(S)$ denotes the class of continuous and bounded functions on $S$. $C_+(S)$ denotes the class of nonnegative continuous functions on $S$, $C_{\text{pol},+} = C_{\text{pol}} \cap C_+$, $C_{\text{b},+} = C_{\text{b}} \cap C_+$ and $C^c_{\text{pol},+} = C_{\text{pol},+} \cap C^c$.

For a function $f : [0, \infty) \to \mathbb{R}^n$, write $|f|^*_t = \sup_{s \in [0,t]} |f(s)|$, and $|f|_t$ for the total variation of $f$ over $[0,t]$ with respect to the Euclidean norm. For a process $X$, we use $X(t)$ and $X_t$ interchangeably.

Denote $B_\varepsilon(x) = \{y \in \mathbb{R}^k : |x - y| < \varepsilon\}$. $S^{d-1}$ denotes the unit sphere in $\mathbb{R}^d$. Infimum over an empty set is regarded as $\infty$. $c, c_1, c_2, \ldots$ denote positive deterministic constants whose values may change from the proof of one result to another.

A function from $[0, \infty)$ to some metric space $E$ is RCLL if it is right-continuous on $[0, \infty)$ and has left limits on $(0, \infty)$. A process is RCLL if, with probability one, its sample paths are RCLL. If $\xi$ is RCLL, denote $\Delta \xi(t) = \xi(t) - \xi(t-)$ for $t > 0$ [see Section 2 for a convention regarding $\Delta \xi(0)$].

**2. Setting and main result.** A filtered probability space $\Phi = (\Omega, \mathcal{F}, (\overline{\mathcal{F}}_t), \mathbb{P}, B)$, satisfying the usual hypotheses, endowed with a $k$-dimensional $(\overline{\mathcal{F}}_t)$-Brownian motion $B$ with drift $b$ and covariance $\Sigma$ is said to be a *system*. Denote by $(\mathcal{F}_t)$ the $\mathbb{P}$-completion of the filtration generated by $B$. We say that $C$ is a cone of $\mathbb{R}^d$ if $C \subset \mathbb{R}^d$, and if $c \in C$ implies $\alpha c \in C$, for all $\alpha \geq 0$. We consider a control problem in which a $p$-dimensional control process $Y$, whose increments take values in a cone $\mathbb{Y}$ (in a sense made precise below), keeps a $k$-dimensional process $X(t) \doteq x + B(t) + GY(t)$ in a cone $\mathbb{X}$, where $G$ is a fixed $k \times p$ matrix of rank $k$ ($k \leq p$). The $k$-dimensional cone $G\mathbb{Y}$ is denoted by $\mathbb{U}$. Our precise assumptions on the cones and related notation are described in what follows. $\mathbb{X}$ (resp., $\mathbb{Y}$, $\mathbb{U}$) is a closed convex cone of $\mathbb{R}^k$ [$\mathbb{R}^p$, $\mathbb{R}^k$] with nonempty interior. It is assumed that

(2.1) $$\mathbb{U} \cap \mathbb{X}^o \neq \varnothing.$$

We remark that, unless $\Sigma$ is degenerate, the above condition is necessary to guarantee the existence of controls; nonetheless, (2.1) will be assumed even for degenerate $\Sigma$. Since $\mathbb{U}$ has nonempty interior, (2.1) implies that there exists a unit vector $\widehat{u}_0 \in \mathbb{U}^o \cap \mathbb{X}^o$. Pick $\widehat{y}_0 \in \mathbb{Y}$ for which $G\widehat{y}_0 = \widehat{u}_0$. The unit vector $\widehat{u}_0$ and the nonzero vector $\widehat{y}_0$ will be fixed throughout. Assume, moreover, that there exist a unit vector $\widehat{u}_1 \in \mathbb{R}^k$, a unit vector $\widehat{y}_1 \in \mathbb{R}^p$, and a constant $a_0 > 0$ such that

$$
\begin{aligned}
u \cdot \widehat{u}_1 \geq a_0 |u|, \quad u \in \mathbb{U}, \quad x \cdot \widehat{u}_1 \geq a_0 |x|, \quad x \in \mathbb{X}, \\
y \cdot \widehat{y}_1 \geq a_0 |y|, \quad y \in \mathbb{Y}.
\end{aligned}
\tag{2.2}
$$

A function $y : [0, \infty) \to \mathbb{R}^p$ is said to have increments in $\mathbb{Y}$ if $y(0) \in \mathbb{Y}$ and $y(t) - y(s) \in \mathbb{Y}$, $0 \leq s < t < \infty$. A process is said to have increments in $\mathbb{Y}$ if, with probability one, its sample paths have increments in $\mathbb{Y}$.

DEFINITION 2.1 (*Admissible control*). An admissible control $Y$ for the system $\Phi$ and the initial data $x \in \mathbb{X}$ is an $(\overline{\mathcal{F}}_t)$-adapted RCLL process with increments in $\mathbb{Y}$ for which the process

$$X(t) \doteq x + B(t) + GY(t), \qquad t \geq 0, \tag{2.3}$$

satisfies $X(t) \in \mathbb{X}$, $t \geq 0$, $\mathbb{P}$-a.s.

By convention, $Y(0-) = 0$ and $X(0-) = x$. The pair $(X, Y)$ (resp., the process $X$) is referred to as an admissible pair (controlled process associated with $Y$) for $\Phi$ and $x$. The class of admissible controls for $\Phi$ and $x$ is denoted by $\bar{\mathcal{A}}(\Phi, x)$, and the class of $(\mathcal{F}_t)$-adapted admissible controls is denoted by $\mathcal{A}(\Phi, x)$. When there is no confusion, we refer to $\bar{\mathcal{A}}(\Phi, x)$ [resp., $\mathcal{A}(\Phi, x)$] as $\bar{\mathcal{A}}(x)$ [$\mathcal{A}(x)$].

Before introducing the cost functional, we fix some notation. Associated with a nondecreasing function $\phi : [0, \infty) \to \mathbb{R}^+$, we define a $\sigma$-finite measure $m$ on $([0, \infty), \mathcal{B}([0, \infty)))$ via the relations $m(a, b] \doteq \phi(b) - \phi(a)$ for $a, b \in (0, \infty)$ and $m\{0\} \doteq \phi(0)$. For $\psi \in L^1([0, \infty), \mathcal{B}([0, \infty)), m)$, we will write $\int_{[0,\infty)} \psi(s) \, d\phi(s) \doteq \int_{[0,\infty)} \psi(s) \, dm(s)$. Note that with this notation, the Stieltjes integral $\int_a^b \psi(s) \, d\phi(s)$ is given as $\int_{(a,b]} \psi(s) \, d\phi(s)$.

Let $y : [0, \infty) \to \mathbb{Y}$ be an RCLL function with increments in $\mathbb{Y}$, and note that, by (2.2), it has bounded variation over finite intervals. Then it can be written as

$$y(t) = \int_{[0,t]} y^\circ(s) \, d|y|(s), \tag{2.4}$$

where $y^\circ$, the Radon–Nikodym derivative $dy/d|y|$, is a measurable function with values in $\mathbb{Y} \cap S^{p-1}$ (see Lemma A.1 in the Appendix). Let $h : \mathbb{Y} \to \mathbb{R}$ be a continuous function satisfying the radial homogeneity condition

$$h(\alpha \eta) = \alpha h(\eta), \qquad \alpha \geq 0, \eta \in \mathbb{Y}. \tag{2.5}$$

Denote

$$\int_{[a,b]} f(s) h(dy(s)) \doteq \int_{[a,b]} f(s) h(y^\circ(s)) \, d|y|(s), \tag{2.6}$$

for every function $f$ for which the right-hand side is well defined. The notation $y^\circ$ and $\int f h(dy)$ of (2.4) and (2.6) is used throughout. In a similar manner, given a system $\Phi$ and an $\{\overline{\mathcal{F}}_t\}$-adapted RCLL process $Y$ with increments in $\mathbb{Y}$, we can find an $\{\overline{\mathcal{F}}_t\}$-progressively measurable process $\{Y_t^\circ\}$ with values in $\mathbb{Y} \cap S^{p-1}$ such that $Y(t) = \int_{[0,t]} Y_s^\circ d|Y|(s)$ (and clearly the value of the integral is independent of the choice of such $\{Y_t^\circ\}$). For this statement, see Lemma A.1. Once more, we will abbreviate $\int_{[a,b]} f(s) h(Y^\circ(s)) d|Y|(s)$ by $\int_{[a,b]} f(s) h(dY(s))$. The cost associated with given system $\Phi$, initial data $x \in \mathbb{X}$ and admissible pair $(X, Y)$ is given as

$$(2.7) \qquad J(\Phi, x, Y) \doteq \mathbb{E} \int_{[0,\infty)} e^{-\beta s} [\ell(X_s)\, ds + h(dY_s)],$$

where, here and throughout, $\mathbb{E}$ denotes expectation with respect to $\mathbb{P}$, and $\beta > 0$ is a constant.

REMARK 2.1. (a) In order to formulate the control problem, one needs to define $h$ only on $\mathbb{Y} \cap S^{p-1}$. However, the radial homogeneous extension of such $h$ to all of $\mathbb{Y}$ will turn out to be convenient.

(b) In the special case where $h$ is linear, say, $h(y) = h_0 \cdot y$, the integral in (2.6) is the same as $\int f(s) h_0 \cdot dy(s)$.

(c) The definition (2.6) reflects the formal identity, $h(dy/d|y|)\, d|y| = h(dy)$, suggested by the radial homogeneity of $h$. Of course, the notation $\int f h(dy)$ should not be confused with a Lebesgue integral against a measure $h$.

Following are our assumptions on $\ell$ and $h$, the function $\ell$ is in $C_+(\mathbb{X})$ and there exist constants $c_{\ell,1}, c_{\ell,2}, c_{\ell,3} \in (0, \infty)$, $m_\ell \in [0, \infty)$ such that

$$(2.8) \qquad c_{\ell,1} |x|^{m_\ell} - c_{\ell,2} \leq \ell(x) \leq c_{\ell,3}(|x|^{m_\ell} + 1), \qquad x \in \mathbb{X}.$$

Note that $m_\ell = 0$ corresponds to the case that $\ell$ is bounded. Let $\mathrm{mod}(r, \delta) \doteq \sup\{\ell(x) - \ell(y) : x, y \in \mathbb{X} \cap B_r(0) : |x - y| \leq \delta\}$ denote the modulus of continuity of $\ell$ on $B_r(0)$. It is assumed that

$$(2.9) \qquad \mathrm{mod}(r+1, \delta) \leq \widehat{m}(\delta)(1 + r^{m_\ell}), \qquad \delta > 0, r > 0,$$

where $\widehat{m}(0+) = 0$. Note that (2.9) is clearly satisfied if $\ell$ is a polynomial. In addition to (2.5), the function $h$ is assumed to be (globally) Lipschitz, convex, and nonnegative on $\mathbb{Y}$. All assumptions mentioned thus far apply throughout this paper. At several places we will also use the following conditions, under which sharper results will be obtained. We will explicitly refer to them when they apply:

$$(2.10) \qquad \text{Either } \Sigma \text{ is nondegenerate or } \widehat{u}_1 \cdot b > 0;$$

$$(2.11) \qquad h(y) \geq c_h |y|, \qquad y \in \mathbb{Y} \text{ for some constant } c_h > 0.$$

We will consider two notions of value function for the control problem. Let

$$\overline{V}(x) \doteq \inf_{\Phi} \inf_{Y \in \overline{\mathcal{A}}(\Phi,x)} J(\Phi, x, Y), \tag{2.12}$$

where in the outer infimum $\Phi$ ranges over all systems, and let

$$V(x) \doteq \inf_{Y \in \mathcal{A}(\Phi,x)} J(\Phi, x, Y). \tag{2.13}$$

Given any two systems $\Phi$ and $\widetilde{\Phi}$, for every $Y \in \mathcal{A}(\Phi, x)$, one can find $\widetilde{Y} \in \mathcal{A}(\widetilde{\Phi}, x)$ such that $(\widetilde{B}, \widetilde{Y})$ is equal in law to $(B, Y)$, and thus $V$ does not depend on $\Phi$.

Consider now the equation

$$((\mathcal{L} + \beta)\psi - \ell) \vee \mathcal{H}(D\psi) = 0, \tag{2.14}$$

where $\mathcal{L}$ denotes the differential operator

$$\mathcal{L} \doteq -\tfrac{1}{2}\text{trace}(\Sigma D^2) - b \cdot D, \tag{2.15}$$

and

$$\mathcal{H}(p) \doteq \sup_{y \in \mathbb{Y}_1} -(Gy \cdot p + h(y)), \qquad p \in \mathbb{R}^k, \tag{2.16}$$

where $\mathbb{Y}_1 = \{y \in \mathbb{Y} : |Gy| = 1\}$.

DEFINITION 2.2 (*Constrained viscosity solution*). (i) For $S \subset \mathbb{X}$, a continuous function $\psi : \overline{S} \to [0, \infty)$ is said to be a viscosity supersolution (resp., subsolution) of (2.14) on $S$ if for all $x \in S$ and all $\varphi \in C^2(S)$ for which $\psi - \varphi$ has a global minimum (maximum) at $x$ one has

$$(\beta\psi(x) + \mathcal{L}\varphi(x) - \ell(x)) \vee \mathcal{H}(D\varphi(x)) \geq 0 \quad (\leq 0).$$

(ii) $\psi$ is said to be a constrained viscosity solution of (2.14) on $\mathbb{X}$ if it is a viscosity supersolution of (2.14) on $\mathbb{X}$ and a viscosity subsolution of (2.14) on $\mathbb{X}^o$.

Our main result characterizes the value function as a constrained viscosity solution of (2.14).

THEOREM 2.1. (i) *If $m_\ell > 0$ (resp., $m_\ell = 0$) then the functions $V$ and $\overline{V}$ are in $C^c_{\text{pol},+}(\mathbb{X})$ [resp., $C_{\text{b},+}(\mathbb{X})$].*

(ii) *Solvability. $V$ and $\overline{V}$ are constrained viscosity solutions of (2.14) on $\mathbb{X}$.*

(iii) *Uniqueness. Let (2.10) hold. Then in the case that $m_\ell > 0$, $V$ is the only such solution in the class $C^c_{\text{pol},+}(\mathbb{X})$; and in case that $m_\ell = 0$, $V$ is the maximal solution in the class $C_{\text{pol},+}(\mathbb{X})$. If, in addition, (2.11) holds, then uniqueness holds in $C_{\text{pol},+}(\mathbb{X})$, $m_\ell \geq 0$.*

(iv) *If* (2.10) *holds, then* $\overline{V} = V$.

REMARK 2.2. Example (a) below shows that condition (2.11) in part (iii) of the theorem is, in a sense, necessary. Example (b) demonstrates the role of the function class $C^c_{\text{pol},+}$. Let $\mathbb{X} = \mathbb{Y} = \mathbb{R}_+$, $G = 1$, $\beta = 1$, $\widehat{u}_1 = 1$ and $h \equiv 0$. In this case the PDE (0.2) reduces to

$$[(\mathcal{L} + 1)f - \ell] \vee (-Df) = 0.$$

(a) *Condition* (2.11) *cannot be dropped in general.* Consider the data $b = 0$, $\Sigma = 1$ and $\ell \equiv 1$. Then (2.10) holds and (2.11) fails. Clearly $m_\ell = 0$ and $V \equiv 1$ on $\mathbb{R}_+$. In this case the statement regarding uniqueness in $C_{\text{pol},+}$ fails to hold, as one checks that, for every $c \in [0,1]$, $\psi \equiv c$ is a constrained viscosity solution of the equation on $\mathbb{R}_+$. This example also demonstrates that in the statement regarding $m_\ell = 0$, one cannot in general replace "maximal solution" by "unique solution."

(b) *The conclusion "uniqueness in $C^c_{\text{pol},+}$" cannot, in general, be replaced by "uniqueness in $C_{\text{pol},+}$."* Assuming now $b = 1$, $\sum = 0$, $l(x) = x + 1$, one finds that $V(x) = x + 2$. Also, it can be easily checked that for every $c \in [0,1]$, $\psi \equiv c$ is a solution. Hence, it follows that, under (2.10) alone, uniqueness in $C_{\text{pol},+}$ does not hold in general.

REMARK 2.3. If $\ell$ is convex and $h$ is linear, the equality of $V$ and $\overline{V}$ is an immediate consequence of Jensen's inequality.

**3. Stochastic networks.** In this section we discuss applications of our result to BCPs. The formulation introduced here is used only in this section, and readers who are not interested in this aspect can safely skip it without losing continuity. As mentioned in the Introduction, BCPs arise from queueing control problems considered in their formal diffusion limit and they often can be transformed into singular control problems of the form studied in this paper. The transformed control problem is sometimes referred to in the literature as the *equivalent workload problem*. Our objective in this section is to describe how our results apply to equivalent workload problems corresponding to a broad family of stochastic networks. To this end, we first define BCPs and quote results of [19] regarding reduction to a singular control problem (no attempt is made to discuss the underlying queueing model or how the BCP arises from it). BCPs were introduced by Harrison in the important work [16] (see [17] for more general formulation). Our presentation follows [19]. Let

$$(3.1) \qquad \widetilde{\Phi} = (\Omega, \mathcal{F}, (\overline{\mathcal{F}}_t), \mathbb{P}, \widetilde{B}),$$

where $(\Omega, \mathcal{F}, (\overline{\mathcal{F}}_t), \mathbb{P})$ is a filtered probability space satisfying the usual hypothesis, and $(\widetilde{B}(t); t \geq 0)$ is an $m$-dimensional $(\overline{\mathcal{F}}_t)$-Brownian motion, with drift $\widetilde{b}$ and covariance $\widetilde{\Sigma}$. The problem data of a BCP is an $m \times n$ matrix $R$, a $p \times n$ matrix $K$ and a vector $z \in \mathbb{R}_+^m$ (termed input–output matrix, capacity consumption matrix and, resp. initial inventory vector). We follow the notation of [19] as far as dimensions of vectors and matrices are concerned, except that we use the symbol $k$ in place of [19]'s $d$. The matrix $K$ is assumed to have rank $p$ ($p \leq n$).

DEFINITION 3.1 (*Admissible control for BCP*). An admissible control $\{L(t); t \geq 0\}$ for the BCP, associated with $\widetilde{\Phi}$ and $z \in \mathbb{R}_+^m$, is an RCLL $(\overline{\mathcal{F}}_t)$-adapted process with values in $\mathbb{R}^n$ such that, setting $Z(t) \doteq z + \widetilde{B}(t) + RL(t)$, $t \geq 0$, and $Y(t) \doteq KL(t)$, one has that $Z(t) \in \mathbb{R}_+^m$ for all $t \geq 0$, and $Y$ has increments in $\mathbb{R}_+^p$.

Denote by $\widetilde{\mathcal{A}}(\widetilde{\Phi}, z)$ the class of all admissible controls for the BCP associated with $\widetilde{\Phi}$ and $z$. The goal is to minimize

$$(3.2) \quad \widetilde{J}(\widetilde{\Phi}, z, L) \doteq \mathbb{E} \int_{[0,\infty)} e^{-\beta t}[\widetilde{\ell}(Z(t))\, dt + h \cdot dY(t)],$$

where $\widetilde{\ell} \in C_+(\mathbb{R}_+^m)$, and $h \in \mathbb{R}_+^p$. Let $\widetilde{V}(z) = \inf_{\widetilde{\Phi}} \inf_{L \in \widetilde{\mathcal{A}}(\widetilde{\Phi}, z)} \widetilde{J}(\widetilde{\Phi}, z, L)$.

Let $\mathcal{B} \doteq \{\lambda \in \mathbb{R}^n : K\lambda = 0\}$, let $\mathcal{N} \doteq R\mathcal{B} \subset \mathbb{R}^m$, and let $q$ be the dimension of $\mathcal{N}$. The dimension of $\mathcal{M} \doteq \mathcal{N}^\perp$ is then $k \doteq m - q$. Let $M$ be any $k \times m$ matrix whose rows span $\mathcal{M}$. By Proposition 2 of [19], there exists a $k \times p$ matrix $G$ such that $MR = GK$. The choice of $G$, in general, is not unique. Set $\mathbb{X} \doteq M\mathbb{R}_+^m$, $\mathbb{Y} = \mathbb{R}_+^p$ and $\mathbb{U} \doteq G\mathbb{Y}$, and note that both $\mathbb{X}$ and $\mathbb{U}$ are subsets of $\mathbb{R}^k$. Define $\ell : \mathbb{X} \to \mathbb{R}_+$ as

$$\ell(x) \doteq \inf\{\widetilde{\ell}(z) : z \in \mathbb{R}_+^m, Mz = x\}.$$

Assume

(3.3) there exists a continuous function $g : \mathbb{X} \to [0, \infty)$ such that $g(x) \in \arg\min\{\widetilde{\ell}(z) : z \in \mathbb{R}_+^m \text{ and } Mz = x\}$, for all $x \in \mathbb{X}$.

With the data $\mathbb{X}$, $\mathbb{Y}$, $G$, $\mathbb{U} = G\mathbb{Y}$, $\Sigma = M\widetilde{\Sigma}M^T$, $b = M\widetilde{b}$, $\ell$, $h$, the singular control problem of Section 2, and in particular, $\overline{\mathcal{A}}(x)$ and $\overline{V}$ are well defined.

THEOREM 3.1 ([19]). *Given $z \in \mathbb{R}_+^m$, let $x = Mz$. Then $\widetilde{V}(z) = \overline{V}(x)$.*

In fact, Harrison and Van Mieghem [19] give an explicit way of constructing an $L$ from a $Y$ such that $\widetilde{J}(\widetilde{\Phi}, z, L) = J(\Phi, x, Y)$ (where $\Phi$ consists of the same filtered probability space as $\widetilde{\Phi}$ and is equipped with the Brownian motion $B$).

We now list some sufficient conditions for our characterization results to hold. For example, if $\widetilde{\ell}$ is linear, nonnegative on $\mathbb{R}_+^m$, and vanishes only at zero, conditions (2.8) and (2.9) are satisfied (in fact, $\ell$ is piecewise linear and $m_\ell = 1$). From Theorem 2 of [3], it follows that (3.3) holds as well. Next, if $G$ and $M$ have full rank, then $\mathbb{X}$ and $\mathbb{U}$ have nonempty interior as subsets of $\mathbb{R}^k$, as required. In, case $M$ and $G$ have nonnegative entries, (2.2) is satisfied with any fixed unit vectors $\widehat{u}_1$ and $\widehat{y}_1$ in the respective positive orthant. Finally, let us assume that $\Sigma$ is nondegenerate. These assumptions hold for a broad family of stochastic networks. Under the *heavy traffic condition* (cf. Assumption 1 of [17]), one can choose $G$ with nonnegative entries (see [17], equation (3.12)). Conditions for a nonnegative choice for $M$ have been given in Theorem 7.3 of [5]. In particular, these conditions hold for open multiclass queueing networks (cf. [5], Section 3.1), parallel server networks (cf. [5], Section 3.2) and several other classes of unitary networks (see [5], Corollary 7.4). As a result, these families of networks are covered under our characterization results.

**4. Preliminary results.** In this section we study some basic properties of the value functions $V$ and $\overline{V}$ [cf. (2.12), (2.13], as well as those of the value function of an analogous problem on a bounded domain, defined below. For $r > 0$, denote

$$(4.1) \qquad \mathbb{X}_r \doteq \{x \in \mathbb{X} : x \cdot \widehat{u}_1 < r\}, \qquad \partial_r \doteq \{x \in \mathbb{X} : x \cdot \widehat{u}_1 = r\}.$$

We will always write $\mathbb{X}_r^c$ for $\mathbb{X} \setminus \mathbb{X}_r$. By (2.2), $\mathbb{X}_r$ and $\partial_r$ are bounded sets. Fix a system $\Phi$, and let $\mathcal{A}(x) = \mathcal{A}(\Phi, x)$. For $Y \in \mathcal{A}(x)$, let $X = x + B + GY$ be the corresponding controlled process, and set $\sigma = \sigma(r)$ as

$$(4.2) \qquad \sigma \doteq \inf\{t : X_t \notin \mathbb{X}_r\}.$$

For $x \in \overline{\mathbb{X}}_r$, let

$$(4.3) \qquad \mathcal{A}_r(x) \doteq \{Y \in \mathcal{A}(x) : \text{on the set } \{\sigma < \infty\}, X_\sigma \in \partial_r\}.$$

Let $\phi$ be any function in $C_+(\partial_r)$. Define for $x \in \overline{\mathbb{X}}_r$ and $Y \in \mathcal{A}_r(x)$ the cost for the bounded domain problem $J_r(x, Y) \equiv J_{r,\phi}(x, Y)$ as

$$(4.4) \quad J_{r,\phi}(x, Y) \doteq \mathbb{E}\left[\int_{[0,\sigma]} e^{-\beta s}[\ell(X_s)\,ds + h(dY_s)] + e^{-\beta \sigma}\phi(X_\sigma)\right],$$

where here and throughout we use the convention that, on the event $\{\sigma = \infty\}$, $[0, \sigma] = [0, \infty)$, and $e^{-\beta \sigma}f(X_\sigma) = 0$. Let also $V_r = V_{r,\phi}$ be defined as

$$(4.5) \qquad V_r(x) \doteq \inf_{Y \in \mathcal{A}_r(x)} J_r(x, Y), \qquad x \in \overline{\mathbb{X}}_r.$$

The notation, $\bar{\mathcal{A}}_r(\Phi, x)$, $J_r(\Phi, x, Y)$ and $\overline{V}_r(x)$ is used analogously to $\bar{\mathcal{A}}(\Phi, x)$, $J(\Phi, x, Y)$ and $\overline{V}(x)$ defined in Section 2.

We first state and prove a result related to the Skorohod problem [11, 12]. For $E = [0, \infty)$ or $E = [0, T]$, let $\mathcal{D}(E : \mathbb{R}^k)$ denote the space of RCLL functions from $E$ to $\mathbb{R}^k$. Denote $\mathcal{D}_{\mathbb{X}}([0, \infty) : \mathbb{R}^k) = \{z \in \mathcal{D}([0, \infty) : \mathbb{R}^k) : z(0) \in \mathbb{X}\}$. Define $\mathcal{D}_{\mathbb{X}}([0, T] : \mathbb{R}^k)$ analogously.

LEMMA 4.1. *There exist maps $\Gamma$ from $\mathcal{D}_{\mathbb{X}}([0, \infty) : \mathbb{R}^k)$ into itself and $\widehat{\Gamma}$ from $\mathcal{D}_{\mathbb{X}}([0, \infty) : \mathbb{R}^k)$ into $\mathcal{D}([0, \infty), \mathbb{R})$ such that the following hold:*

(i) *If $z \in \mathcal{D}_{\mathbb{X}}([0, \infty) : \mathbb{R}^k)$, $v = \widehat{\Gamma}(z)$ and $x = \Gamma(z)$ then $v : [0, \infty) \to \mathbb{R}_+$ is RCLL and nondecreasing, $x = z + \widehat{u}_0 v = z + G\widehat{y}_0 v$, and $x(t) \in \mathbb{X}$ for all $t \geq 0$.*

(ii) *The maps $\widehat{\Gamma}$ and $\Gamma$ are Lipschitz in the following sense. For all $z_1, z_2 \in \mathcal{D}_{\mathbb{X}}([0, \infty) : \mathbb{R}^k)$ and all $T > 0$,*

$$|\widehat{\Gamma}(z_1) - \widehat{\Gamma}(z_2)|_T^* + |\Gamma(z_1) - \Gamma(z_2)|_T^* \leq \kappa |z_1 - z_2|_T^*,$$

*where the constant $\kappa < \infty$ does not depend on $T$, $z_1$ and $z_2$.*

(iii) *The map $\widehat{\Gamma}$ is nonanticipating in the following sense. For every $T \in (0, \infty)$ there exists a map $\widehat{\Gamma}_T : \mathcal{D}_{\mathbb{X}}([0, T] : \mathbb{R}^k) \to \mathcal{D}([0, T] : \mathbb{R})$, such that for $z \in \mathcal{D}_{\mathbb{X}}([0, \infty) : \mathbb{R}^k)$,*

$$\widehat{\Gamma}(z)|_{[0,T]} = \widehat{\Gamma}_T(z|_{[0,T]}).$$

(iv) *If $z \in \mathcal{D}_{\mathbb{X}}([0, \infty) : \mathbb{R}^k)$ takes values in $\mathbb{X}$ on $[0, T]$, then $v = \widehat{\Gamma}(z)$ vanishes on $[0, T]$.*

The notation $\Gamma$ and $\widehat{\Gamma}$ is kept throughout this paper.

REMARK 4.1. Clearly, $\Gamma(z) = z + G\widehat{y}_0 \widehat{\Gamma}(z)$. The lemma will help us construct one admissible control from another, as follows. If $X = x + B + GY$ is the controlled process for a $Y \in \mathcal{A}(x)$, then letting $\zeta = \widehat{\Gamma}(\widetilde{x} + B + GY)$, the process $\widetilde{Y} = Y + \widehat{y}_0 \zeta$ is seen to be admissible for $\widetilde{x}$.

PROOF OF LEMMA 4.1. Let $\overline{S}$ be the collection of all unit vectors $\overline{s} \in \mathbb{R}^k$ such that $\{\xi \in \mathbb{R}^k : \overline{s} \cdot \xi \geq 0\} \supset \mathbb{X}$. Then one has the following representation for $\mathbb{X}$ (cf. [30]):

$$\mathbb{X} = \bigcap_{\overline{s} \in \overline{S}} \{\xi \in \mathbb{R}^k : \overline{s} \cdot \xi \geq 0\}.$$

Recall that $\widehat{u}_0 \in \mathbb{X}^o$. For $\xi \in \mathbb{R}^k$, let $\pi(\xi)$ denote the projection of $\xi$ onto the boundary $\partial \mathbb{X}$ along $\widehat{u}_0$. Explicitly, $\pi(\xi) = \xi + \alpha(\xi)\widehat{u}_0$, where

$$\alpha(\xi) = \sup_{\overline{s} \in \overline{S}} -\xi \cdot \overline{s} / \widehat{u}_0 \cdot \overline{s}.$$

It is elementary to check that the range of $\pi$ is $\partial \mathbb{X}$ and that $\pi$ is globally Lipschitz. Given $z \in \mathcal{D}_{\mathbb{X}}([0,\infty):\mathbb{R}^k)$, let

$$v(t) = 0 \vee \sup_{0 \leq s \leq t} \widehat{u}_0 \cdot (\pi(z(s)) - z(s)), \qquad t \geq 0,$$

and $x = z + \widehat{u}_0 v$. Then using the fact that $\pi(z(t)) = \pi(x(t))$, it is not hard to check that $\widehat{u}_0 \cdot (x(t) - \pi(x(t))) \geq 0$, $t \geq 0$. Thus, $x(t) = \pi(x(t)) - \alpha(x(t))\widehat{u}_0$, where $\alpha(x(t)) \leq 0$. Since both $\widehat{u}_0$ and $\pi(x(t))$ are in $\mathbb{X}$, this shows that $x(t) \in \mathbb{X}$, and part (i) of the lemma is established. Parts (ii), (iii) and (iv) follow by construction and the fact [used in proof of (ii)] that $\pi$ is globally Lipschitz. $\square$

Recall that $h$ is radially homogeneous and convex. This is easily seen to imply that

(4.6) $$h(y+z) \leq h(y) + h(z), \qquad y, z \in \mathbb{Y}.$$

As a result we have the following.

LEMMA 4.2. *Let $Y_1$ and $Y_2$ have increments in $\mathbb{Y}$ and set $Y = Y_1 + Y_2$. Then for all $t \in (0, \infty)$,*

(4.7) $$\int_{[0,t]} e^{-\beta s} h(dY(s)) \leq \int_{[0,t]} e^{-\beta s} h(dY_1(s)) + \int_{[0,t]} e^{-\beta s} h(dY_2(s)).$$

*Let, for $i=1,2$, $\mu_i$ be a $\sigma$-finite measure on $(\mathbb{R}_+, \mathcal{B}(\mathbb{R}_+))$ defined as*

$$\mu_i(B) \doteq \int_B d|Y_i|(s), \qquad B \in \mathcal{B}(\mathbb{R}_+),$$

*where $\mathcal{B}(\mathbb{R}_+)$ is the Borel $\sigma$-field on $\mathbb{R}_+$. If $\mu_1$ and $\mu_2$ are mutually singular, then* (4.7) *holds with equality.*

PROOF. Along the lines of [15], page 320, let $\widehat{Y}(s) = \widehat{y}_1 \cdot Y(s)$ and denote $d\xi(s) = d|Y|(s)$ and $d\widehat{\xi}(s) = d\widehat{Y}(s)$. Also let

$$\mu(B) = \int_B d\xi, \qquad \widehat{\mu}(B) = \int_B d\widehat{\xi},$$

for Borel $B \subset [0,\infty)$, and define similarly $\widehat{Y}_i$, $d\xi_i$, $d\widehat{\xi}_i$ and $\widehat{\mu}_i$ for $i=1,2$. By (2.2), $a_0|y| \leq \widehat{y}_1 \cdot y \leq |y|$ for $y \in \mathbb{Y}$. Hence, $\mu$ and $\widehat{\mu}$ are mutually absolutely continuous, and a similar statement holds for $\mu_i$ and $\widehat{\mu}_i$, $i=1,2$. Since $\widehat{Y} = \widehat{Y}_1 + \widehat{Y}_2$ and $\widehat{Y}_i$ are nondecreasing, clearly, $\widehat{\mu}_i$ is absolutely continuous with respect to $\widehat{\mu}$, for $i=1,2$. This shows that $\mu_i$ is absolutely continuous with respect to $\mu$, and we denote by $d\xi_i/d\xi$ the respective Radon–Nikodym derivatives. Thus,

$$Y^\circ = Y_1^\circ \frac{d\xi_1}{d\xi} + Y_2^\circ \frac{d\xi_2}{d\xi}, \qquad \mu\text{-a.e.}$$

By (4.6),

$$h(Y^\circ(s)) \leq \sum_i \frac{d\xi_i}{d\xi}(s) h(Y_i^\circ(s)), \qquad \mu\text{-a.e.} \tag{4.8}$$

and as a result

$$\int_{[0,t]} e^{-\beta s} h(Y^\circ(s)) \, d\xi(s) \leq \sum_i \int_{[0,t]} e^{-\beta s} h(Y_i^\circ(s)) \, d\xi_i(s). \tag{4.9}$$

Next, suppose that $\mu_1$ and $\mu_2$ are mutually singular. Then there are disjoint sets $S_i \in \mathcal{B}(\mathbb{R}_+)$ such that $\mu_i(S_i^c) = 0$. Clearly, $d\xi_i/d\xi = \mathbb{1}_{S_i} d\xi_i/d\xi, \mu$-a.e.. Since (4.6) holds with equality whenever either $y$ or $z$ vanishes, (4.8) holds with equality, and so does (4.9). □

LEMMA 4.3. *Let $x \in \mathbb{X}$ and $y \in \mathbb{Y}$. If $x + Gy \in \mathbb{X}$, then*

$$V(x + Gy) + h(y) \geq V(x), \qquad \overline{V}(x + Gy) + h(y) \geq \overline{V}(x).$$

PROOF. We will only prove the result for $V$. The proof for $\overline{V}$ is similar. Assume without loss of generality that $\mathcal{A}(x + Gy)$ is nonempty. Let $Y \in \mathcal{A}(x + Gy)$. Then the corresponding controlled process $X = x + Gy + B + GY$ takes values in $\mathbb{X}$. Set $\overline{Y} = y + Y$ and $\overline{X} = x + B + G\overline{Y}$. Clearly, $\overline{X} = X$, and $\overline{Y}$ has increments in $\mathbb{Y}$ and, therefore, $\overline{Y} \in \mathcal{A}(x)$. Also, $\Delta \overline{Y}(0) = y + \Delta Y(0)$, and therefore, by (4.6), $h(\Delta \overline{Y}(0)) \leq h(y) + h(\Delta Y(0))$. Thus,

$$J(x + Gy, Y) \geq J(x, \overline{Y}) - h(y) \geq V(x) - h(y).$$

Since $Y \in \mathcal{A}(x + Gy)$ is arbitrary, the result follows. □

In the proof of the next result, and several times in the paper, we use the fact that if $X$ is the controlled process corresponding to some $Y \in \mathcal{A}(x)$, then, by (2.2), for $0 \leq t \leq s < \infty$,

$$(4.10) \quad |X_s| \geq \widehat{u}_1 \cdot X_s \geq \widehat{u}_1 \cdot X_t + \widehat{u}_1 \cdot (B_s - B_t) \geq a_0 |X_t| - |B_s - B_t|.$$

As immediate consequences of this inequality, we have

$$(4.11) \quad 2^m \mathbb{E}|X_s|^m \geq a_0^m \mathbb{E}|X_t|^m - 2^m \mathbb{E}|B_t - B_s|^m, \qquad s \geq t, \ m \geq 0,$$

$$(4.12) \qquad |X|_s^* \leq 2a_0^{-1}(|X_s| + |B|_s^*), \qquad s \geq 0.$$

LEMMA 4.4. *There exist constants $a_1, a_2, a_3 > 0$ such that*

$$a_1 |x|^{m_\ell} - a_2 \leq V(x) \leq a_3(1 + |x|^{m_\ell}), \qquad x \in \mathbb{X}.$$

*The above inequality also holds with $V$ replaced by $\overline{V}$.*

PROOF. Once more we will only prove the result for $V$. Let $x \in \mathbb{X}$ and define $Y = \widehat{y}_0 \widehat{\Gamma}(x+B)$ and $X = \Gamma(x+B)$. Then $X = x + B + GY$, and it follows from parts (i) and (iii) of Lemma 4.1 that $Y \in \mathcal{A}(x)$.

Denoting by $\mathbf{0}$ the zero trajectory in $\mathbb{R}^k$, it is clear that $\Gamma(x+\mathbf{0})(t) = x$, $t \geq 0$. Hence, by Lemma 4.1(ii), $|X(t) - x| + |Y(t)| \leq c_1 |B|_t^*$. With the notation of the proof of Lemma 4.2, the Radon–Nikodym derivative of $\mu$ w.r.t. $\widehat{\mu}$ is bounded above by $a_0^{-1}$, as follows from (2.2). Thus,

$$\int_{[0,\infty)} e^{-\beta s} d|Y|(s) \leq a_0^{-1} \int_{[0,\infty)} e^{-\beta s} d\widehat{Y}(s)$$

(4.13)
$$= (a_0 \beta)^{-1} \int_{[0,\infty)} e^{-\beta s} \widehat{Y}(s)\, ds$$

$$\leq (a_0 \beta)^{-1} \int_{[0,\infty)} e^{-\beta s} |Y(s)|\, ds.$$

By (2.8) and (4.13),

$$J(x,Y) = \mathbb{E}\int_{[0,\infty)} e^{-\beta t}(\ell(X(t))\, dt + h(Y^\circ(t))\, d|Y|(t))$$

$$\leq c_2 \mathbb{E}\int_{[0,\infty)} e^{-\beta t}(1 + |X(t)|^{m_\ell} + |Y|(t))\, dt$$

$$\leq c_3 \mathbb{E}\int_{[0,\infty)} e^{-\beta t}(1 + |x|^{m_\ell} + (|B|_t^*)^{m_\ell \vee 1})\, dt$$

$$\leq c_4(1 + |x|^{m_\ell}),$$

where $c_2, c_3, c_4$ do not depend on $x$. By (2.8) and (4.11), for every admissible $Y$,

$$J(x,Y) \geq c_5 \mathbb{E}\int_{[0,\infty)} e^{-\beta t} |X_t|^{m_\ell}\, dt - c_6 \geq c_7 |x|^{m_\ell} - c_8,$$

where $c_7, c_8 > 0$ are independent of $x$ and $Y$. $\square$

LEMMA 4.5. *$V$ and $\overline{V}$ are continuous on $\mathbb{X}$.*

PROOF. We will only consider $V$. Fix $r > 0$, and given arbitrary $\varepsilon \in (0,1)$, consider $x_1, x_2 \in \mathbb{X} \cap B_r(0)$ with $|x_1 - x_2| < \delta < 1$, where $\delta > 0$ will be chosen later. Fix $Y_1 \in \mathcal{A}(x_1)$ such that

(4.14) $$J(x_1, Y_1) \leq V(x_1) + \varepsilon/2,$$

and let $X_1$ be the corresponding controlled process, namely, $X_1 = x_1 + B + GY_1$. Let $Z = \widehat{\Gamma}(x_2 + B + GY_1)$, $Y_2 = Y_1 + \widehat{y}_0 Z$ and $X_2 = x_2 + B + GY_2 =$ $x_2 + B + GY_1 + \widehat{u}_0 Z$. By Lemma 4.1, $Y_2 \in \mathcal{A}(x_2)$ and $X_2$ is the corresponding controlled process. Note that

$$X_1 = \Gamma(X_1), \qquad X_2 = \Gamma(x_2 - x_1 + X_1).$$

Hence, by Lemma 4.1,

(4.15) $$|X_1(t) - X_2(t)| + |Y_1(t) - Y_2(t)| \leq c_1 \delta,$$

where $c_1$ does not depend on $t$, $x_1$ and $x_2$. Assume that $\delta$ is small enough so that $c_1 \delta < 1$. In particular, this shows that $\int_{[0,\infty)} e^{-\beta t} \, dZ(t) \leq c_2 \delta$, and thus, by (2.8), (2.9), Lemmas 4.2 and 4.4, we have

$$J(x_2, Y_2) - J(x_1, Y_1) \leq \mathbb{E} \int_{[0,\infty)} e^{-\beta t} [(\ell(X_2(t)) - \ell(X_1(t))) \, dt + h(\widehat{y}_0) \, dZ(t)]$$

$$\leq \mathbb{E} \int_{[0,\infty)} e^{-\beta t} \operatorname{mod}(|X_1(t)| + 1, \delta) \, dt + c_3 \delta$$

(4.16) $$\leq c_4 \widehat{m}(c_1 \delta) \mathbb{E} \int_{[0,\infty)} e^{-\beta t} (\ell(X_1(t)) + 1) \, dt + c_3 \delta$$

$$\leq c_5 \widehat{m}(c_1 \delta)(V(x_1) + 1) + c_3 \delta$$

$$\leq c_6 \widehat{m}(c_1 \delta)(r^{m_\ell} + 1) + c_3 \delta.$$

Choosing $\delta$ so small that the expression on the last line is bounded by $\varepsilon/2$, we conclude that $V(x_2) - V(x_1) \leq \varepsilon$ whenever $|x_1 - x_2| \leq \delta$ and $x_1, x_2 \in \mathbb{X} \cap B(0, r)$. Since $r$ is arbitrary, $V$ is continuous on $\mathbb{X}$. $\square$

REMARK 4.2. Lemmas 4.4 and 4.5 are seen to imply part (i) of Theorem 2.1.

For $x \in \mathbb{X}_r$ and $v \in \mathbb{U}$ for which there exists $\rho > 0$ such that $x + \rho v \in \mathbb{X}_r^c$, let

(4.17) $$\gamma_r(x, v) \doteq \inf\{\rho > 0 : x + \rho v \in \mathbb{X}_r^c\}.$$

Clearly, $x + \gamma_r(x, v) v \in \partial_r$ for $x, v$ as above, and

(4.18) $$x \in \mathbb{X}_r, v \in \mathbb{U}, \qquad x + v \in \mathbb{X}_r^c \implies \gamma_r(x, v) \leq 1.$$

LEMMA 4.6. *Let $\phi \in C_+(\partial_r)$ and suppose that (2.10) holds. Then $V_{r,\phi}$ is continuous on $\overline{\mathbb{X}}_r$.*

PROOF. Below we use implicitly the fact that, in the definition of $V_r$, the values $X$ and $Y$ take on the interval $(\sigma, \infty)$ are immaterial. Fix $\varepsilon > 0$ and consider all $x_1, x_2 \in \overline{\mathbb{X}}_r$ with $|x_1 - x_2| \leq \delta \leq 1/2$, where $\delta > 0$ will be chosen later. Let $Y_1 \in \mathcal{A}_r(x_1)$ be such that

(4.19) $$J_{r,\phi}(x_1, Y_1) \leq V_{r,\phi}(x_1) + \varepsilon/2,$$

let $X_1$ be the corresponding controlled process, $X_1 \equiv X^{x_1} \doteq x_1 + B + GY_1$, and let $\sigma_1$ be the corresponding exit time. Let $Z = \widehat{\Gamma}(x_2 + B + GY_1)$ and $X = \Gamma(x_2 + B + GY_1)$. By Lemma 4.1(i), (iii), $Z$ is $\mathbb{R}_+$-valued, adapted, RCLL and nondecreasing, and $X(t) \in \mathbb{X}$, $t \geq 0$. Also, by Lemma 4.1(iv), $\widehat{\Gamma}(x_1 + B + GY_1) = 0$, hence, by Lemma 4.1(ii),

$$(4.20) \qquad 0 \leq Z(t) \leq \kappa\delta, \qquad t \geq 0.$$

Let $Y = Y_1 + \widehat{y}_0 Z$. Define $\tau \doteq \inf\{t : X(t) \in \mathbb{X}_r^c\}$. Let $\sigma_2 = \sigma_1 \wedge \tau$ and set $(X_2, Y_2) = (X, Y)$ on $[0, \sigma_2)$. Define also

$$(4.21) \quad Y_2(\sigma_2) \doteq \begin{cases} Y_2(\tau-) + \gamma_r(X(\tau-), \Delta X(\tau))\Delta Y(\tau), & \tau \leq \sigma_1, \\ Y(\sigma_1) + \gamma_r(X(\sigma_1), \widehat{u}_0)\widehat{y}_0, & \sigma_1 < \tau, \end{cases}$$

and

$$(4.22) \qquad X_2(\sigma_2) = X_2(\sigma_2-) + G\Delta Y_2(\sigma_2).$$

Note that $X_2(\sigma_2) \in \partial_r$. Hence, $X_2$ is well defined until $\sigma_2$, the first time it exits $\mathbb{X}_r$. We leave $X_2$ undefined on $(\sigma_2, \infty)$. Below we sometimes write $U_i$ for $GY_i$, $i = 1, 2$.

Note that $\sigma_2 \leq \sigma_1$. By (4.17) and (4.21), the random variable $\Delta X_2(\sigma_2)$ is $\mathbb{U}$-valued and $\mathcal{F}_{\sigma_2}$-measurable, and since

$$X_2(t) = x_2 + B(t) + GY_2(t), \qquad t \leq \sigma_2,$$

we see that $Y_2 \in \mathcal{A}(x_2)$ and $X_2$ is the corresponding controlled process. On the time interval $[0, \sigma_2)$ we have

$$X_1 = \Gamma(X_1), \qquad X_2 = \Gamma(x_2 - x_1 + X_1),$$

and therefore, by Lemma 4.1,

$$(4.23) \quad |X_1(t) - X_2(t)| + |Y_1(t) - Y_2(t)| \leq c_1|x_1 - x_2| \leq c_1\delta, \qquad t < \sigma_2,$$

and, similarly,

$$(4.24) \qquad |X_1(t) - X(t)| \leq c_1\delta, \qquad t \leq \sigma_2,$$

where $c_1$ does not depend on $t, x_1, x_2$ and $\delta$.

We show that there exists a constant $c$ not depending on $x_1, x_2$ and $\delta$ such that

$$(4.25) \qquad |X_1(\sigma_2) - X_2(\sigma_2)| + h(\Delta Y_2(\sigma_2)) - h(\Delta Y_1(\sigma_2)) \leq c\delta.$$

Different arguments are used in the two cases below:

*Case* (i): $\tau \leq \sigma_1$. By (4.17) and (4.21), $X(\tau) - X_2(\tau) \in \mathbb{U}$. Moreover, since $X_2(\tau) \in \partial_r$ and $X_1(\tau) \in \overline{\mathbb{X}}_r$, $\widehat{u}_1 \cdot (X_2(\tau) - X_1(\tau)) \geq 0$. Therefore, by (2.2) and (4.24),

$$
\begin{aligned}
a_0 |X(\tau) - X_2(\tau)| &\leq \widehat{u}_1 \cdot (X(\tau) - X_2(\tau)) \\
&\leq \widehat{u}_1 \cdot (X(\tau) - X_1(\tau)) \leq c_1 \delta.
\end{aligned}
\tag{4.26}
$$

Combining (4.24) and (4.26), we get

$$
|X_1(\sigma_2) - X_2(\sigma_2)| \leq c_2 \delta. \tag{4.27}
$$

Since $X(\tau-) \in \mathbb{X}_r$ and $X(\tau) \in \mathbb{X}_r^c$, it follows from (4.18) that $\gamma \doteq \gamma_r(X(\tau-), \Delta X(\tau)) \leq 1$. With (4.21), we have $\Delta Y_2(\sigma_2) = \gamma \Delta Y(\sigma_2) = \gamma \Delta Y_1(\sigma_2) + \gamma \widehat{y}_0 \Delta Z(\sigma_2)$. Thus, by the Lipschitz property of $h$ and (4.20),

$$
h(\Delta Y_2(\sigma_2)) \leq h(\gamma \Delta Y_1(\sigma_2)) + c_3 \gamma |\widehat{y}_0| \Delta Z(\sigma_2) \leq h(\Delta Y_1(\sigma_2)) + c_4 \delta, \tag{4.28}
$$

and (4.25) follows.

*Case* (ii): $\sigma_1 < \tau$. Note that $\sigma_1 = \sigma_2$ and by (4.24), $|X_1(\sigma_2) - X(\sigma_2)| \leq c_1 \delta$. Since $X_1(\sigma_2), X_2(\sigma_2) \in \partial_r$, we have $\widehat{u}_1 \cdot X_1(\sigma_2) = \widehat{u}_1 \cdot X_2(\sigma_2)$. Also, note that $X_2(\sigma_2) - X(\sigma_2) \in \mathbb{X}$. Thus, using (4.24):

$$
\begin{aligned}
a_0 |X_2(\sigma_2) - X(\sigma_2)| &\leq \widehat{u}_1 \cdot (X_2(\sigma_2) - X(\sigma_2)) \\
&= \widehat{u}_1 \cdot (X_1(\sigma_2) - X(\sigma_2)) \leq c_1 \delta, \qquad \sigma_1 < \tau.
\end{aligned}
\tag{4.29}
$$

Combining (4.24) and (4.29), the estimate (4.27) holds for an appropriate constant $c_2$. Now, by (4.21) and (4.22),

$$
\Delta X_2(\sigma_1) = G \Delta Y_2(\sigma_1) = G(Y_1(\sigma_1-) - Y_2(\sigma_1-)) + G \Delta Y_1(\sigma_1) + \gamma' \widehat{u}_0,
$$

where $\gamma' = \gamma_r(X(\sigma_1), \widehat{u}_0) + Z(\sigma_1)$. Since $\Delta X_1(\sigma_1) = G \Delta Y_1(\sigma_1)$, we have

$$
\gamma' \leq |\Delta X_2(\sigma_1) - \Delta X_1(\sigma_1)| + c_5 |Y_1(\sigma_1-) - Y_2(\sigma_1-)|,
$$

and since $\sigma_1 = \sigma_2$, we have, by (4.23) and (4.27) that

$$
\gamma' \leq c_6 \delta. \tag{4.30}
$$

Using again (4.21), $\Delta Y_2(\sigma_2) - \Delta Y_1(\sigma_2) = \gamma' \widehat{y}_0$. Combining (4.30) and the Lipschitz property of $h$, we establish (4.25).

Let

$$
\begin{aligned}
m(\alpha) &= \max\{|\ell(y) - \ell(z)| : y, z \in \overline{\mathbb{X}}_r, |y - z| \leq \alpha\} \\
&\vee \max\{|\phi(y) - \phi(z)| : y, z \in \partial_r, |y - z| \leq \alpha\}.
\end{aligned}
$$

Recalling the convention $e^{-\beta\sigma}f(X_\sigma) = 0$ when $\sigma = \infty$,

$$J_{r,\phi}(x_2, Y_2) - J_{r,\phi}(x_1, Y_1)$$

$$\leq \mathbb{E}\int_{[0,\sigma_2]} e^{-\beta t}[(\ell(X_2(t)) - \ell(X_1(t)))\,dt$$

$$+ h(dY_2(t)) - h(dY_1(t))]$$

(4.31)
$$+ \mathbb{E}\{\mathbb{1}_{\sigma_1 = \sigma_2} e^{-\beta\sigma_2}(\phi(X_2(\sigma_2)) - \phi(X_1(\sigma_2)))\}$$

$$+ \mathbb{E}\{\mathbb{1}_{\sigma_2 < \sigma_1}[e^{-\beta\sigma_2}\phi(X_2(\sigma_2)) - e^{-\beta\sigma_1}\phi(X_1(\sigma_1))]\}$$

$$\leq c_7\delta + c_8 m(c_8\delta)$$

$$+ \mathbb{E}\{\mathbb{1}_{\sigma_2 < \sigma_1}[e^{-\beta\sigma_2}\phi(X_2(\sigma_2)) - e^{-\beta\sigma_1}\phi(X_1(\sigma_1))]\}.$$

The estimate $c_7\delta$ in the last line above follows on using Lemma 4.2, (4.20) and (4.25). By (4.27), since $X_2(\sigma_2) \in \partial_r$, we have $r - \widehat{u}_1 \cdot X_1(\sigma_2) \leq c_9\delta$. Thus, for $\alpha > 0$, $\sigma_1 > \sigma_2 + \alpha$ implies that, for $t \in [\sigma_2, \sigma_2 + \alpha]$,

$$r > \widehat{u}_1 \cdot X_1(t) \geq \widehat{u}_1 \cdot x_1 + \widehat{u}_1 \cdot B(t) + \widehat{u}_1 \cdot U_1(\sigma_2)$$

$$= \widehat{u}_1 \cdot X_1(\sigma_2) + \widehat{u}_1 \cdot (B(t) - B(\sigma_2))$$

$$\geq r - c_9\delta + \widehat{u}_1 \cdot (B(t) - B(\sigma_2)).$$

Therefore,

(4.32)
$$P(\sigma_1 > \sigma_2 + \alpha) \leq P\left(\max_{0 \leq t \leq \alpha} \widehat{u}_1 \cdot (B(\sigma_2 + t) - B(\sigma_2)) \leq c_9\delta\right)$$

$$\doteq \lambda_1(\alpha, \delta).$$

Note that, if either $\Sigma$ is nondegenerate or $\widehat{u}_1 \cdot b > 0$, we have that, for each fixed $\alpha$, $\lambda_1(\alpha, \delta) \to 0$ as $\delta \to 0$. Moreover, writing $X_1(\sigma_1) - X_1(\sigma_2) = B(\sigma_1) - B(\sigma_2) + U_1(\sigma_1) - U_1(\sigma_2)$, we have $\widehat{u}_1 \cdot (U_1(\sigma_1) - U_1(\sigma_2)) \leq c_9\delta + |B(\sigma_1) - B(\sigma_2)|$. Using (2.2), on $\{0 < \sigma_1 - \sigma_2 < \alpha\}$,

(4.33)
$$|X_1(\sigma_1) - X_1(\sigma_2)| \leq c_{10}\delta + c_{10}\max\{|B(\sigma_2 + t) - B(\sigma_2)| : t \in [0, \alpha]\}$$

$$= c_{10}\delta + \lambda_2(\alpha).$$

Using (4.32) and (4.33) in (4.31), we obtain

$$J_{r,\phi}(x_2, U_2) - J_{r,\phi}(x_1, U_1)$$

$$\leq c_{11}(\delta + m(c_4\delta) + \lambda_1(\alpha, \delta) + \alpha + \mathbb{E}(m(c_{10}\delta + \lambda_2(\alpha))))$$

$$\leq c_{11}(\delta + m(c_4\delta) + \lambda_1(\alpha, \delta) + \alpha + m(2c_{10}\delta) + \mathbb{E}(m(2\lambda_2(\alpha)))).$$

Note that $\mathbb{E}(m(2\lambda_2(\alpha))) \to 0$ as $\alpha \to 0$ since $m_2$ is uniformly bounded (bound only depends on $r$) and $\lambda_2(\alpha) \to 0$ as $\alpha \to 0$. Choose $\alpha$ small enough and then $\delta$ small enough so that the right-hand side is bounded by $\varepsilon/2$. Combining this with (4.19), we have that $V_{r,\phi}(x_2) \leq V_{r,\phi}(x_1) + \varepsilon$ whenever $|x_1 - x_2| \leq \delta$. This proves the continuity of $V_{r,\phi}$ on $\overline{\mathbb{X}}_r$. □

The following lemma shows that the infimum in the definition of $V(x)$ [$\overline{V}(x)$] in (2.13) [resp. (2.12)] can equivalently be performed over a class of admissible controls under which the controlled process's moments are finite and sub-exponential in the time variable.

LEMMA 4.7. *For $x \in \mathbb{X}$ and a system $\Phi$, let*

$$\mathcal{A}_F(x) = \Big\{Y \in \mathcal{A}(x) : \forall\ \alpha, t > 0\ \mathbb{E}|X_Y(t)|^\alpha < \infty;$$

$$\forall\ \alpha > 0 \lim_{t \to \infty} e^{-\beta t} \mathbb{E}|X_Y(t)|^\alpha = 0\Big\},$$

*where $X_Y$ is the controlled process corresponding to $Y$. Define $\bar{\mathcal{A}}_F(\Phi, x)$ similarly by replacing $\mathcal{A}(\Phi, x)$ above by $\bar{\mathcal{A}}(\Phi, x)$. Then*

(4.34) $$V(x) = \inf\{J(x, \Phi, Y) : \in \mathcal{A}_F(\Phi, x)\}$$

*and*

(4.35) $$\overline{V}(x) = \inf_\Phi \inf\{J(x, \Phi, Y) : Y \in \bar{\mathcal{A}}_F(\Phi, x)\}.$$

PROOF. We will prove (4.34). The proof of (4.35) is identical. Fix $x \in \mathbb{X}$. Given $p > |x|$, an admissible control $Y$ and the corresponding controlled process $X$ for which $J(x, Y) < \infty$, let $\lambda_p = \inf\{t : \widehat{u}_1 \cdot X(t) \geq p\}$. Note that $\lambda_p$ may assume the value 0. Let

$$Z_p(t) = X(t)\mathbb{1}_{[0,\lambda_p)}(t) + (X(\lambda_p-) + B(t) - B(\lambda_p))\mathbb{1}_{[\lambda_p,\infty)}(t), \qquad t \geq 0,$$

where $\mathbb{1}_{[0,\lambda_p)}$ is understood as zero in case $\lambda_p = 0$ and $X(0-)$, by convention, equals $x$. Set

(4.36)
$$Y_p = Y\mathbb{1}_{[0,\lambda_p)} + \widehat{y}_0 \widehat{\Gamma}(Z_p),$$
$$X_p \doteq \Gamma(Z_p) = x + B + GY\mathbb{1}_{[0,\lambda_p)} + \widehat{u}_0 \widehat{\Gamma}(Z_p).$$

By Lemma 4.1, $Y_p$ is an admissible control and $X_p$ is the corresponding controlled process. Denote $\widetilde{X}_p(t) = X(t)\mathbb{1}_{[0,\lambda_p)}(t) + X(\lambda_p-)\mathbb{1}_{[\lambda_p,\infty)}$. Clearly, $\widetilde{X}_p$ takes values in $\mathbb{X}$ and, therefore, $\widetilde{X}_p = \Gamma(\widetilde{X}_p)$. Thus, by Lemma 4.1, on $\{\lambda_p < \infty\}$,

(4.37)
$$|X_p - \widetilde{X}_p|^*_T \leq \kappa |Z_p - \widetilde{X}_p|^*_T$$
$$= \kappa \sup\{|B(s) - B(\lambda_p)| : s \in [\lambda_p, T]\}, \qquad T \geq \lambda_p.$$

Thus, by (4.36), the right-hand side of (4.37) is an upper bound for $|X_p(T) - X_p(\lambda_p-)|$ for $T \geq \lambda_p$. A similar argument using the fact that $\widehat{\Gamma}(\widetilde{X}_p) = 0$ shows that the right-hand side of (4.37) is an upper bound also for $|Y_p(T) - Y_p(\lambda_p-)|$. Thus, have on $\{\lambda_p < \infty\}$,

$$
(4.38) \quad \begin{aligned} &|X_p(T) - X_p(\lambda_p-)| + |Y_p(T) - Y_p(\lambda_p-)| \\ &\leq c_1 \sup\{|B(t) - B(\lambda_p)| : t \in [\lambda_p, T]\}, \quad T \geq \lambda_p. \end{aligned}
$$

Since on $\{\lambda_p < \infty\}$, we have $|X_p(\lambda_p-)| = |X(\lambda_p-)| \leq a_0^{-1}p$, the sub-exponential behavior of the $\alpha$th moment of $X_p(t)$ follows from a similar property of the Brownian motion. Thus, to prove (4.34) it suffices to show that

$$
(4.39) \quad \limsup_{p \to \infty} J(x, Y_p) \leq J(x, Y), \quad Y \in \mathcal{A}(x).
$$

Note that, by (2.8),

$$
\mathbb{1}_{[\lambda_p, \infty)}(s) c_{\ell,1} |X(s)|^{m_\ell} \leq \ell(X_s) + c_{\ell,2}
$$

and since $J(x, Y) < \infty$, we have $\mathbb{E} \int_0^\infty e^{-\beta s}(\ell(X_s) + c_{\ell,2})\,ds < \infty$. Also, clearly, $\lambda_p \to \infty$, a.s., as $p \to \infty$. Thus, $\int_{\lambda_p}^\infty e^{-\beta s}|X(s)|^{m_\ell}\,ds \to 0$ a.s. and in $L^1$. Using (4.10) with $t = \lambda_p-$, this shows that

$$
(4.40) \quad \lim_{p \to \infty} \mathbb{E}[e^{-\beta \lambda_p}|X(\lambda_p-)|^{m_\ell}] = 0.
$$

Note that, with $\xi \doteq \widehat{y}_0 \widehat{\Gamma}(Z_p)$, $[Y_p(t) - Y_p(\lambda_p-)]\mathbb{1}_{t \geq \lambda_p} = \xi(t)$. Hence, using (4.38), we get

$$
\begin{aligned} \int_{[\lambda_p, \infty)} e^{-\beta s} h(d\xi(s)) &\leq c \int_{[\lambda_p, \infty)} e^{-\beta s}\,d|\xi|(s) \\ &\leq c \int_{[\lambda_p, \infty)} e^{-\beta s} |Y_p(t) - Y_p(\lambda_p-)|\,ds. \end{aligned}
$$

Therefore, with (2.8) and (4.38), we have

$$
\begin{aligned} &J(x, Y_p) - J(x, Y) \\ &\leq c_{\ell,3} \mathbb{E} \int_{[\lambda_p, \infty)} e^{-\beta s}[(|X_p(s)|^{m_\ell} + 1)\,ds + h(d\xi(s))] \\ &\leq c_2 \mathbb{E} \int_{[\lambda_p, \infty)} e^{-\beta s}\left(|X_p(\lambda_p-)|^{m_\ell} + \sup_{t \in [\lambda_p, s]} |B_s - B_{\lambda_p}|^{m_\ell \vee 1} + 1\right) ds, \end{aligned}
$$

where the constant $c_2$ does not depend on $p$ and $Y$. Using the above along with (4.40) yields (4.39). This proves the lemma. □

The following lemma shows that, in computing the value function, the class of admissible controls can be further restricted.

LEMMA 4.8. *Fix $x \in \mathbb{X}$, a system $\Phi$ and $Y \in \bar{\mathcal{A}}_F(\Phi, x)$ with $J(x, \Phi, Y) < \infty$. Then one can find a $c \in (0, \infty)$ such that, for all $\varepsilon \in (0, 1)$, setting*

$$T^\varepsilon = \inf\left\{t : \mathbb{E}\int_{[t,\infty)} e^{-\beta s}(\ell(X_s)\,ds + h(dY_s)) \le \varepsilon\right\} \vee \log(\varepsilon^{-1/\beta}), \quad (4.41)$$

$$Z^\varepsilon(t) = X(t)\mathbb{1}_{[0,T^\varepsilon]}(t) + (X(T^\varepsilon) + B(t) - B(T^\varepsilon))\mathbb{1}_{(T^\varepsilon,\infty)}(t), \quad (4.42)$$

$$t \ge 0,$$

*and*

$$Y^\varepsilon = Y\mathbb{1}_{[0,T^\varepsilon]} + \widehat{y}_0\widehat{\Gamma}(Z^\varepsilon),$$
$$X^\varepsilon = \Gamma(Z^\varepsilon) = x + B + GY\mathbb{1}_{[0,T^\varepsilon]} + \widehat{u}_0\widehat{\Gamma}(Z^\varepsilon), \quad (4.43)$$

*we have that, for all $T \ge T^\varepsilon$,*

$$(4.44) \quad |X_T^\varepsilon - X_{T^\varepsilon}| + |Y_T^\varepsilon - Y_{T^\varepsilon}^\varepsilon| \le c\sup\{|B(t) - B(T^\varepsilon)| : t \in [T^\varepsilon, T]\}$$

*and $J(x, Y^\varepsilon) \le J(x, Y) + c\varepsilon$.*

PROOF. Note that $X^\varepsilon$ is the controlled process corresponding to $Y^\varepsilon$. Since $X^\varepsilon(t) \in \mathbb{X}$ for $t \le T^\varepsilon$, we have $\widehat{\Gamma}(Z^\varepsilon)(t) = 0$ for $t \le T^\varepsilon$ and so

$$(4.45) \qquad (Y^\varepsilon(t), X^\varepsilon(t)) = (Y(t), X(t)), \qquad t \le T^\varepsilon.$$

Arguing again by Lemma 4.1, letting $\widetilde{X}^\varepsilon = X(\cdot \wedge T^\varepsilon)$, we have $\widetilde{X}^\varepsilon = \Gamma(\widetilde{X}^\varepsilon)$ and therefore, for $T \ge T^\varepsilon$,

$$|X^\varepsilon - \widetilde{X}^\varepsilon|_T^* \le \kappa|Z^\varepsilon - \widetilde{X}^\varepsilon|_T^* = \kappa\sup\{|B(t) - B(T^\varepsilon)| : t \in [T^\varepsilon, T]\}.$$

Combining the above with (4.43) yields (4.44).

Next, moment estimates on the Brownian motion imply

$$(4.46) \quad \mathbb{E}|X^\varepsilon(T)|^{m_\ell} \le c_1(\mathbb{E}|X_{T^\varepsilon}|^{m_\ell} + (T - T^\varepsilon)^m + 1), \qquad T \ge T^\varepsilon,$$

where $c_1, m$ do not depend on $T, \varepsilon$ and $Y$. Hence,

$$(4.47) \qquad \mathbb{E}\int_{(T^\varepsilon,\infty)} e^{-\beta s}|X^\varepsilon(s)|^{m_\ell}\,ds \le c_2 e^{-\beta T^\varepsilon}(\mathbb{E}|X_{T^\varepsilon}|^{m_\ell} + 1).$$

Arguing as in the previous lemma, by (4.44), we also have

$$\mathbb{E}\int_{[T^\varepsilon,\infty)} e^{-\beta s}h(dY_s^\varepsilon) \le c\mathbb{E}\int_{[T^\varepsilon,\infty)} e^{-\beta s}|Y_s^\varepsilon - Y_{T^\varepsilon}^\varepsilon|\,ds$$
$$(4.48) \qquad \qquad \le c_3 e^{-\beta T^\varepsilon} \le c_3\varepsilon.$$

Next, by (4.41), (4.11) and (2.8),

$$\varepsilon \geq \mathbb{E}\int_{[T^\varepsilon, T^\varepsilon+1)} e^{-\beta s}(\ell(X_s)\,ds + h(dY_s))$$

$$\geq e^{-\beta(T^\varepsilon+1)}\mathbb{E}\int_{T^\varepsilon}^{T^\varepsilon+1}(c_{\ell,1}|X_s|^{m_\ell} - c_{\ell,2})\,ds$$

$$\geq c_4 e^{-\beta T^\varepsilon}[\mathbb{E}|X_{T^\varepsilon}|^{m_\ell} - c_5(1+\mathbb{E}(|B|_1^*)^{m_\ell})]$$

$$= c_4 e^{-\beta T^\varepsilon}[\mathbb{E}|X_{T^\varepsilon}|^{m_\ell} - c_6],$$

where $c_4, c_5, c_6 > 0$ are independent of $\varepsilon$ and $Y$. Hence, by (4.41),

(4.49) $$e^{-\beta T^\varepsilon}\mathbb{E}|X_{T^\varepsilon}|^{m_\ell} \leq \frac{\varepsilon}{c_4} + \frac{c_6}{c_4}e^{-\beta T^\varepsilon} \leq c_7\varepsilon.$$

By (4.45), (4.47), (4.48) and (4.49),

$$J(x, Y^\varepsilon) \leq \mathbb{E}\int_{[0,T^\varepsilon]} e^{-\beta s}(\ell(X_s^\varepsilon)\,ds + h(dY_s^\varepsilon)) + c_8 e^{-\beta T^\varepsilon}(\mathbb{E}|X_{T^\varepsilon}|^{m_\ell} + 1) + c_3\varepsilon$$

$$\leq J(x, Y) + c_9\varepsilon.$$

This proves the lemma. □

**5. Solvability.** In this section we establish part (ii) of Theorem 2.1 by proving the following result.

THEOREM 5.1. *Both $V$ and $\overline{V}$ are constrained viscosity solutions of (2.14) on $\mathbb{X}$.*

A key to the proof will be the following DPPs, Propositions 5.1 and 5.2, whose proof, along with the proof of Lemma 5.1, are postponed to Section 8. The first DPP regards $V$.

PROPOSITION 5.1. *Let $x \in \mathbb{X}$ and $\varepsilon > 0$ be given. For $Y \in \mathcal{A}(x)$ and $X$ the controlled process associated with $x$ and $Y$, let*

(5.1) $$\tau_Y \equiv \tau \doteq \inf\{t \geq 0 : X(t) \notin B_\varepsilon(x)\}.$$

*Then for $t \in [0, \infty)$,*

(5.2) $$V(x) = \inf_{Y \in \mathcal{A}(x)} \mathbb{E}\bigg[\int_{[0, t \wedge \tau]} e^{-\beta s}(\ell(X(s))\,ds + h(dY(s))) + e^{-\beta(t \wedge \tau)}V(X(t \wedge \tau))\bigg].$$

In order to present the DPP associated with $\overline{V}$, we need to introduce some notation. Let $\Phi = (\Omega, \mathcal{F}, (\overline{\mathcal{F}}_t), \mathbb{P}, B)$ be a system and $\zeta$ be an $\mathbb{R}^k$ valued random variable with probability distribution $\mu$ given on $(\Omega, \overline{\mathcal{F}}_0, \mathbb{P})$. We will refer to $(\Phi, \zeta)$ as an *extended system*. An $\overline{\mathcal{F}}_t$ adapted RCLL process $Y$ with increments in $\mathbb{Y}$ is said to be an admissible control for $(\Phi, \zeta)$ if $X_t \doteq \zeta + B_t + GY_t$ satisfies $X_t \in \mathbb{X}$, $\mathbb{P}$ a.s., for all $t \geq 0$. Denote the class of all such admissible controls by $\bar{\mathcal{A}}(\Phi, \zeta)$. Let

$$(5.3) \qquad J(\Phi, \zeta, Y) \doteq \mathbb{E} \int_{[0,\infty)} e^{-\beta s} [\ell(X_s)\, ds + h(dY_s)].$$

Denote the class of all $\overline{\mathcal{F}}_t$ adapted RCLL processes with increments in $\mathbb{Y}$ by $\mathcal{A}_0(\Phi)$. Given $Y^0 \in \mathcal{A}_0(\Phi)$ and $\zeta$ as above, let

$$(5.4) \qquad Z_t \doteq \zeta + B_t + GY^0_t, \qquad \eta \doteq \widehat{\Gamma}(Z) \quad \text{and} \quad Y \doteq Y^0 + \widehat{y}_0 \eta.$$

Then, $Y \in \bar{\mathcal{A}}(\Phi, \zeta)$ and the corresponding controlled process $X$ is given as $X_t = \zeta + B_t + GY^0_t + G\widehat{y}_0 \eta_t$. Let

$$(5.5) \qquad \widehat{\mathcal{F}}_t \doteq \sigma\{\zeta, X_s, \eta_s : 0 \leq s \leq t\}.$$

Let

$$(5.6) \qquad \begin{aligned}\widehat{\mathcal{A}}(\Phi, \zeta) = \{Y \in \bar{\mathcal{A}}(\Phi, \zeta) | \text{ there exists } Y^0 \in \mathcal{A}_0(\Phi) \text{ satisfying (5.4)} \\ \text{and } Y^0 \text{ is } \widehat{\mathcal{F}}_t \text{ adapted}\}.\end{aligned}$$

Elements of $\widehat{\mathcal{A}}(\Phi, \zeta)$ will be referred to as feedback controls. When $\mu = \delta_x$, for some $x \in \mathbb{X}$, we will write $\widehat{\mathcal{A}}(\Phi, \zeta)$ as $\widehat{\mathcal{A}}(\Phi, x)$. Let, for $x \in \mathbb{X}$,

$$(5.7) \qquad \widehat{V}(x) \doteq \inf_{\Phi} \inf_{Y \in \widehat{\mathcal{A}}(\Phi,x)} J(x, \Phi, Y).$$

The following lemma proves the equality of $\widehat{V}$ and $\overline{V}$.

LEMMA 5.1. *For all $x \in \mathbb{X}, \overline{V}(x) = \widehat{V}(x)$.*

We can now state the second dynamic programming principle.

PROPOSITION 5.2. *Let $x \in \mathbb{X}$ and $\varepsilon > 0$ be given. For $Y \in \bar{\mathcal{A}}(\Phi, x)$, let $\tau_Y$ be defined via (5.1). Then for $t \in (0, \infty)$,*

$$(5.8) \qquad \begin{aligned}\overline{V}(x) = \inf_{\Phi} \inf_{Y \in \widehat{\mathcal{A}}(\Phi,x)} \mathbb{E}\bigg[&\int_{[0,t\wedge\tau]} e^{-\beta s}(\ell(X(s))\, ds + h(dY(s))) \\ &+ e^{-\beta(t\wedge\tau)} \overline{V}(X(t\wedge\tau))\bigg].\end{aligned}$$

REMARK 5.1. When using Propositions 5.1 and 5.2 we can assume without loss of generality, that the infimum is taken only on those $Y \in \mathcal{A}(x)$ [resp.

$\widehat{\mathcal{A}}(\Phi, x)]$ for which, on the set $\{\tau \leq t\}$, $X(\tau) \in \partial B_\varepsilon(x)$. More precisely, in case of Proposition 5.1, for $t \in (0, \infty)$, let

$$\mathcal{A}_{1,t}(x) \doteq \{Y \in \mathcal{A}(x) \colon \text{On the set } \{\tau \leq t\}, X(\tau) \in \partial B_\varepsilon(x)\}.$$

Then

(5.2) $$V(x) = \inf_{Y \in \mathcal{A}_{1,t}(x)} \mathbb{E}\left[\int_{[0, t \wedge \tau]} e^{-\beta s}(\ell(X(s))\,ds + h(dY(s))) + e^{-\beta(t \wedge \tau)} V(X(t \wedge \tau))\right].$$

Indeed, for a general $Y \in \mathcal{A}(x)$, consider $(\widetilde{X}, \widetilde{Y})$ that agree with $(X, Y)$ on $[0, \tau)$ and, on $\{\tau \leq t\}$, satisfy $\Delta \widetilde{Y}(\tau) = \alpha \Delta Y(\tau)$, where $\alpha \in (0, 1]$ is such that $\widetilde{X}(\tau) \in \partial B_\varepsilon(x)$. Denoting by $J_1(Y)$ the expectation on the right-hand side of (5.2) and writing $\delta \doteq \Delta Y(\tau) - \Delta \widetilde{Y}(\tau)$, clearly, $J_1(Y) - J_1(\widetilde{Y}) = E\{\mathbb{1}_{\tau \leq t} e^{-\tau}[h(\delta) + V(\widetilde{X}(\tau) + G\delta) - V(\widetilde{X}(\tau))]\} \geq 0$ by Lemma 4.3. A similar statement holds for the dynamic programming principle in Proposition 5.2.

The proofs of the two dynamic programming principles and Lemma 5.1 are deferred to Section 8.

PROPOSITION 5.3. *Both $V$ and $\overline{V}$ are viscosity supersolutions of* (2.14) *on* $\mathbb{X}$.

PROOF. We will first consider $V$. Fix $x \in \mathbb{X}$ and let $\varphi \in C^2(\mathbb{X})$ be such that $V - \varphi$ has a global minimum at $x$. We can assume, without loss of generality, that $V(x) - \varphi(x) = 0$. We need to show that either

(5.9) $$\beta \varphi(x) + \mathcal{L}\varphi(x) - \ell(x) \geq 0$$

or

(5.10) $$\inf\{Gy \cdot D\varphi(x) + h(y) \colon y \in \mathbb{Y}_1\} \leq 0.$$

Arguing by contradiction, assume that neither of the above assertions is true. Then one can find $\theta > 0$ and $\varepsilon > 0$ such that, for all $\bar{x} \in \overline{B_\varepsilon(x)} \cap \mathbb{X}$,

(5.11) $$\beta \varphi(\bar{x}) + \mathcal{L}\varphi(\bar{x}) - \ell(\bar{x}) \leq -\theta$$

and $Gy \cdot D\varphi(\bar{x}) + h(y) \geq \theta$ for all $y \in \mathbb{Y}_1$. The latter implies that, for all $y \in \mathbb{Y}$,

(5.12) $$Gy \cdot D\varphi(\bar{x}) + h(y) \geq \theta|Gy|.$$

Let $t > 0$, fix $Y \in \mathcal{A}_{1,t}(x)$ and denote $U = GY$. Let $X$ be the corresponding controlled process. Denote

$$Y^{\mathrm{c}}(t) = Y(t) - \sum_{0 \leq s \leq t} \Delta Y(s), \qquad U^{\mathrm{c}} = GY^{\mathrm{c}}.$$

Then $U^c$ is continuous, $U^c(0) = 0$, and it has increments in $\mathbb{U}$. Let $\tau$ be as in (5.1). An application of Itô's formula gives

$$\varphi(x) = \mathbb{E}[e^{-\beta(t\wedge\tau)}\varphi(X_{t\wedge\tau})] + \mathbb{E}\int_0^{t\wedge\tau} e^{-\beta s}(\mathcal{L}\varphi(X_s) + \beta\varphi(X_s))\,ds$$

(5.13)
$$- \mathbb{E}\int_{[0,t\wedge\tau]} e^{-\beta s} D\varphi(X_{s-}) \cdot dU_s^c$$

$$- \mathbb{E}\sum_{0\le s\le t\wedge\tau} e^{-\beta s}(\varphi(X_s) - \varphi(X_{s-})).$$

From (5.12), it follows that, for $0 \le s \le t \wedge \tau$,

$$\varphi(X_s) - \varphi(X_{s-}) = \int_0^1 D\varphi(X_{s-} + \sigma\Delta U_s) \cdot G\Delta Y_s\,d\sigma$$

(5.14)
$$\ge \theta|\Delta U_s| - h(\Delta Y_s).$$

We also have, by (5.12),

$$\int_{[0,t\wedge\tau]} e^{-\beta s} D\varphi(X_{s-}) \cdot dU_s^c$$

(5.15)
$$\ge \theta \int_{[0,t\wedge\tau]} e^{-\beta s}\,d|U^c|_s - \int_{[0,t\wedge\tau]} e^{-\beta s} h(dY_s^c).$$

Combining (5.11), (5.14), (5.14), (5.15) and the second part of Lemma 4.2, we obtain

$$\varphi(x) \le \mathbb{E}[e^{-\beta(t\wedge\tau)}\varphi(X_{t\wedge\tau})] + \mathbb{E}\int_{[0,t\wedge\tau]} e^{-\beta s}[\ell(X_s)\,ds + h(dY_s)]$$

$$- \theta\mathbb{E}\int_{[0,t\wedge\tau]} e^{-\beta s}(ds + d|U|_s).$$

Since $\varphi \le V$ and $|U_r| \le \int_{[0,r]} d|U|_s$, we have

$$\varphi(x) \le \mathbb{E}[e^{-\beta(t\wedge\tau)}V(X_{t\wedge\tau})] + \mathbb{E}\int_{[0,t\wedge\tau]} e^{-\beta s}[\ell(X_s)\,ds + h(dY_s)]$$

(5.16)
$$- \theta e^{-\beta t}\mathbb{E}[t \wedge \tau + |U_{t\wedge\tau}|].$$

Taking infimum over all $Y \in \mathcal{A}_{1,t}(x)$, we have from Proposition 5.1 and Remark 5.1 that

$$V(x) = \varphi(x) \le V(x) - \theta e^{-\beta t}\alpha(t),$$

where

$$\alpha(t) \doteq \inf_{Y \in \mathcal{A}_{1,t}(x)} \mathbb{E}[t \wedge \tau + |U_{t\wedge\tau}|].$$

A contradiction will be obtained by showing

(5.17) $$\exists t > 0 \quad \text{s.t.} \quad \alpha(t) > 0.$$

Recall that $X = x + U + B$, $U = GY$ and that, since $Y \in \mathcal{A}_{1,t}(x)$, one has $|X_{t \wedge \tau} - x| = \varepsilon$ on $\{\tau \leq t\}$. Hence,

$$\mathbb{E}|U_{t \wedge \tau}|\mathbb{1}_{\tau \leq t} \geq \mathbb{E}[(\varepsilon - |B_{t \wedge \tau}|)\mathbb{1}_{\tau \leq t}].$$

Thus,

$$\mathbb{E}[t \wedge \tau + |U_{t \wedge \tau}|] \geq \mathbb{E}[t\mathbb{1}_{\tau > t} + (\varepsilon/2)\mathbb{1}_{\tau \leq t, |B|_t^* < \varepsilon/2}] \geq [t \wedge (\varepsilon/2)]P(|B|_t^* < \varepsilon/2).$$

Clearly, for all $t > 0$ small enough, $P(|B|_t^* < \varepsilon/2) > 0$. Note that the qualifier "small enough" is needed since the Brownian motion is allowed to be degenerate. This proves (5.17) and hence the first part of the result. The proof that $\overline{V}$ is a supersolution as well is similar; instead of taking infimum over all $Y \in \mathcal{A}_{1,t}(x)$ in (5.16), we take infimum over all feedback controls and use the dynamic programming principle in Proposition 5.2. $\square$

PROPOSITION 5.4. *Both $V$ and $\overline{V}$ are viscosity subsolutions of* (2.14) *on $\mathbb{X}^o$.*

PROOF. Once more we only prove the result for $V$. Fix $x \in \mathbb{X}^o$ and let $\varphi \in C^2(\mathbb{X}^o)$ be such that $V - \varphi$ has a global maximum at $x$. We need to show that

(5.18) $$\beta V(x) + \mathcal{L}\varphi(x) - \ell(x) \leq 0$$

and

(5.19) $$Gy \cdot D\varphi(x) + h(y) \geq 0, \qquad y \in \mathbb{Y}_1.$$

We can assume without loss of generality that $\varphi(x) = V(x)$. Thus, $V \leq \varphi$ on $\mathbb{X}^o$. For all $\delta > 0$ small enough, one has $x + \delta Gy \in \mathbb{X}^o$ for all $y \in \mathbb{Y}_1$. Hence, by Lemma 4.3,

$$\varphi(x + \delta Gy) - \varphi(x) \geq V(x + \delta Gy) - V(x) \geq -\delta h(y), \qquad y \in \mathbb{Y}_1.$$

Dividing by $\delta$ and taking $\delta \to 0$ proves (5.19).

To prove (5.18), let $\varepsilon > 0$ be such that $\overline{B}_\varepsilon(x) \subset S^o$. For a control $Y$ and a corresponding controlled process $X$, let

$$\tau_Y^\varepsilon \doteq 1 \wedge \inf\{t : X(t) \notin \overline{B}_\varepsilon(x)\}.$$

Now set $Y_s = 0$ for $s \in [0, \tau_Y^\varepsilon]$ and, thus, $X_s = x + B_s$ for $s \leq \tau_Y^\varepsilon$. In what follows denote $\tau = \tau_Y^\varepsilon$. An application of Itô's formula gives

(5.20) $$V(x) = \varphi(x) = \mathbb{E}[e^{-\beta \tau}\varphi(X_\tau)] + \mathbb{E}\int_0^\tau e^{-\beta s}(\mathcal{L}\varphi(X_s) + \beta\varphi(X_s))\,ds.$$

Using Proposition 5.1, the inequality $V \leq \varphi$ and (5.20),

$$V(x) \leq \mathbb{E}\left(\int_0^\tau e^{-\beta s} \ell(X_s)\,ds + e^{-\beta \tau}\varphi(X_\tau)\right)$$
(5.21)
$$= \varphi(x) + \mathbb{E}\int_0^\tau e^{-\beta s}(\ell(X_s) - \mathcal{L}\varphi(X_s) - \beta\varphi(X_s))\,ds.$$

Recalling that $V(x) = \varphi(x)$ and denoting

$$\zeta(\bar{x}) \doteq \ell(\bar{x}) - \mathcal{L}\varphi(\bar{x}) - \beta\varphi(\bar{x}), \qquad \bar{x} \in \overline{B}_\varepsilon(x),$$

we have $\mathbb{E}\int_0^\tau e^{-\beta s}\zeta(X_s)\,ds \geq 0$. Hence,

$$\zeta(x)\mathbb{E}\int_0^\tau e^{-\beta s}\,ds \geq -\mathbb{E}\int_0^\tau e^{-\beta s}(\zeta(X_s) - \zeta(x))\,ds \geq -\alpha(x,\varepsilon)\mathbb{E}\int_0^\tau e^{-\beta s}\,ds,$$

where

$$\alpha(x,\varepsilon) = \max_{\bar{x}\in\overline{B}_\varepsilon(x)} |\zeta(\bar{x}) - \zeta(x)|.$$

Since $\tau > 0$ a.s., it follows that $\zeta(x) \geq -\alpha(x,\varepsilon)$. Taking $\varepsilon \to 0$, we obtain $\zeta(x) \geq 0$, proving (5.18) and hence the first part of the result. The second part once more is obtained upon using Proposition 5.2 rather than Proposition 5.1 in proving the statement analogous to (5.22). □

Combining Propositions 5.3 and 5.4, we obtain Theorem 5.1.

**6. Uniqueness on bounded domain.** Recall the notation $\mathbb{X}_r, \partial_r, \sigma = \sigma(r)$, as well as $J_{r,\phi}, V_{r,\phi}$ and $\overline{V}_{r,\phi}$ from Section 4. In particular, recall that

$$V_{r,\phi}(x) = \inf_{Y\in\mathcal{A}_r(x)} J_{r,\phi}(x,Y), \qquad \overline{V}_{r,\phi}(x) = \inf_{Y\in\overline{\mathcal{A}}_r(x)} J_{r,\phi}(x,Y).$$

Let $r > 0$ be fixed throughout this section. In this section we prove the following.

PROPOSITION 6.1. *Let $u \in C_+(\mathbb{X})$ be a constrained viscosity solution of (2.14) on $\mathbb{X}$. Assume that (2.10) holds. Then for $x \in \overline{\mathbb{X}}_r$,*

(6.1)                $u(x) = V_{r,u}(x) = \overline{V}_{r,u}(x).$

Recall Definition 2.2 of viscosity sub- and supersolutions.

DEFINITION 6.1. Let $\phi:\partial_r \to \mathbb{R}$ be given. We say that a continuous function $\psi:\overline{\mathbb{X}}_r \to \mathbb{R}_+$ is a constrained viscosity solution of (2.14) on $\mathbb{X}_r$ with the Dirichlet boundary condition $\phi$ if the following conditions hold:
  (i) $\psi$ is a viscosity supersolution of (2.14) on $\mathbb{X}_r$;
  (ii) $\psi$ is a viscosity subsolution of (2.14) on $\mathbb{X}_r^o$;

(iii) $\psi = \phi$ on $\partial_r$.

REMARK 6.1. If $u : \mathbb{X} \to \mathbb{R}_+$ is a constrained viscosity solution of (2.14) on $\mathbb{X}$, then it is clearly a constrained viscosity solution of (2.14) on $\mathbb{X}_r$ with the Dirichlet boundary condition $u|_{\partial_r}$.

The proof of the following result is similar to that of Theorem 5.1 and therefore is omitted.

THEOREM 6.1. *Let $\phi \in C_+(\partial_r)$. Assume that (2.10) holds. Then both $V_{r,\phi}$ and $\overline{V}_{r,\phi}$ are constrained viscosity solutions of (2.14) on $\mathbb{X}_r$ with the Dirichlet boundary condition $\phi$.*

The following will be the principal tool in proving Proposition 6.1.

THEOREM 6.2. *Let $\phi \in C_+(\partial_r)$. Let $u, v \in C_+(\overline{\mathbb{X}}_r)$ be two constrained viscosity solutions of (2.14) on $\mathbb{X}_r$ with the Dirichlet boundary condition $\phi$. Then $u = v$.*

PROOF OF PROPOSITION 6.1. The function $V_{r,u}$ is clearly nonnegative, and by Lemma 4.6 and Theorem 6.1, it is a $C_+(\overline{\mathbb{X}}_r)$ constrained viscosity solution of (2.14) on $\mathbb{X}_r$ with the Dirichlet boundary condition $u|_{\partial_r}$. In view of Remark 6.1, so is $u$. A similar statement holds for $\overline{V}_{r,u}$. Hence, the result follows from Theorem 6.2. □

In the rest of this section we prove Theorem 6.2.

LEMMA 6.1. *For every $\xi \in \overline{\mathbb{X}}_r$, there exist $\eta = \eta(\xi) \in \mathbb{R}^n$ and $a = a(\xi) > 0$ such that*

(6.2) $\qquad B_{ta}(x + t\eta) \subset \mathbb{X}_r^o \qquad \forall\, x \in \overline{\mathbb{X}}_r \cap B_a(\xi) \text{ and } \forall\, t \in (0, 1].$

PROOF. Fix $\xi \in \overline{\mathbb{X}}_r$. Let $w \in \mathbb{X}_r^o$ and let $a > 0$ be so small that $B_{2a}(w) \subset \mathbb{X}_r^o$. If $x \in \overline{\mathbb{X}}_r$ and $y \in \mathbb{X}_r^o$, then convexity of $\mathbb{X}_r$ implies that any point on the line segment joining $x$ and $y$, excluding $x$, belongs to the interior $\mathbb{X}_r^o$. Hence, for any $x \in \overline{\mathbb{X}}_r$,

$$\bigcup_{t \in (0,1]} B_{2ta}(x + t(w - x)) = \{x + t(z - x) : z \in B_{2a}(w), t \in (0, 1]\} \subset \mathbb{X}_r^o.$$

The inclusion above may be written as

(6.3) $$\bigcup_{t \in (0,1]} B_{2ta}(x + t(w - x)) \subset \mathbb{X}_r^o.$$

Let $\eta = w - \xi$. Clearly, for $x \in B_a(\xi) \cap \overline{\mathbb{X}}_r$, $\bigcup_{t \in (0,1]} B_{ta}(x+t\eta) \subset \bigcup_{t \in (0,1]} B_{2ta}(x+t(w-x))$. Hence, the result follows from (6.3). $\square$

PROOF OF THEOREM 6.2. We introduce some notation specific to the proof. For $n \in \mathbb{N}$, denote by $\mathcal{S}(n)$ the space of symmetric $n \times n$ matrices. We write $M_1 \leq M_2$ if and only if $M_2 - M_1$ is nonnegative definite, $M_1, M_2 \in \mathcal{S}(n)$. Let $S$ be relatively open in $\mathbb{X}$. For $x \in \mathbb{X}$, $r \in \mathbb{R}$, $p \in \mathbb{R}^k$ and $A \in \mathcal{S}(k)$, write

$$F(x, r, p, A) = \beta r - b \cdot p - \tfrac{1}{2}\mathrm{trace}(\Sigma A) - \ell(x).$$

For $x \in S$ and a real valued continuous function $\psi$ on $\overline{\mathbb{X}}_r$, denote

$$J_S^{2,+}\psi(x) \doteq \{(D\varphi(x), D^2\varphi(x)) : \varphi \in C^2(S) \text{ and}$$
$$\psi - \varphi \text{ has a local maximum at } x\},$$

$$J_S^{2,-}\psi(x) \doteq \{(D\varphi(x), D^2\varphi(x)) : \varphi \in C^2(S) \text{ and}$$
$$\psi - \varphi \text{ has a local minimum at } x\}.$$

Define the closures of the above sets in the following way. For $x \in S$,

$$\overline{J}_S^{2,+}\psi(x) \doteq \{(p, M) \in \mathbb{R}^k \times \mathcal{S}(k) : \text{ there exists a sequence}$$
$$(x_n, p_n, M_n) \in \mathbb{X}_r \times \mathbb{R}^k \times \mathcal{S}(k) \text{ s.t.}$$
$$(p_n, M_n) \in J_S^{2,+}\psi(x_n), \text{ and}$$
$$(x_n, \psi(x_n), p_n, M_n) \to (x, \psi(x), p, M)$$
$$\text{as } n \to \infty\}.$$

Define $\overline{J}_S^{2,-}\psi(x)$ analogously. For short, write $J^{2,+}$ for $J_{\mathbb{X}_r^\circ}^{2,+}$ and similarly define $J^{2,-}$, $\overline{J}^{2,+}$ and $\overline{J}^{2,-}$. Generic elements of $\mathbb{Y}$ will be denoted by $\underline{y}$ (rather than $y$).

To prove the result, it suffices to show that

$$(6.4) \qquad u(x) \leq v(x), \qquad x \in \overline{\mathbb{X}}_r.$$

Let $c_1 \doteq a_0^{-1}$. Note that $c_1 G\underline{y} \cdot \widehat{u}_1 \geq 1$, $\underline{y} \in \mathbb{Y}_1$ [cf. (2.2)]. Let $\psi(x) \doteq c_1 x \cdot \widehat{u}_1 - c_2$, $x \in \overline{\mathbb{X}}_r$, where $c_2$ is fixed and large enough so that

$$(6.5) \quad \psi(x) \leq -1 \quad \text{and} \quad (\mathcal{L} + \beta)\psi(x) - \ell(x) \leq -1, \qquad x \in \overline{\mathbb{X}}_r.$$

Thus,

$$(6.6) \qquad G\underline{y} \cdot D\psi(x) \geq 1, \qquad x \in \overline{\mathbb{X}}_r, \underline{y} \in \mathbb{Y}_1.$$

For $\alpha \in (0, 1)$, define

$$(6.7) \qquad u_\alpha \doteq \alpha u + (1-\alpha)\psi.$$

It suffices to show that, for every $\alpha \in (0,1)$,

(6.8) $$u_\alpha(x) \leq v(x), \qquad x \in \overline{\mathbb{X}}_r.$$

We argue by contradiction and assume that (6.8) does not hold. Therefore, there exist $\alpha \in (0,1)$ and $\xi \in \overline{\mathbb{X}}_r$ such that

(6.9) $$u_\alpha(\xi) - v(\xi) = \max_{x \in \overline{\mathbb{X}}_r}(u_\alpha(x) - v(x)) \doteq \delta > 0.$$

We first argue that if $z \in \mathbb{X}_r^o$ and $(p, A) \in \overline{J}^{2,+} u_\alpha(z)$, then

(6.10) $$F(z, u_\alpha(z), p, A) \leq -(1-\alpha), \qquad \mathcal{H}(p) \leq -(1-\alpha).$$

To this end, consider first $(p, A) \in J^{2,+} u_\alpha(z)$. Let $p^*$ and $A^*$ be such that

(6.11) $$\begin{aligned} p &= \alpha p^* + (1-\alpha) D\psi(z), \\ A &= \alpha A^* + (1-\alpha) D^2 \psi(z) = (1-\alpha) A^*. \end{aligned}$$

Then by (6.7), $(p^*, A^*) \in J^{2,+} u(z)$. Using the subsolution property of $u$ (cf. Definition 6.1), we have

(6.12) $$F(z, u(z), p^*, A^*) \leq 0.$$

Noting that the map $(r, q, X) \mapsto F(z, r, q, X)$ is affine and combining (6.5), (6.7), (6.11) and (6.12), we obtain the first inequality in (6.10). Now since $z \in \mathbb{X}_r^o$ and $(p, A) \in J^{2,+} u_\alpha(z)$ are arbitrary, and $F$ is continuous in all variables, the first inequality in (6.10) holds, in fact, for all $z \in \mathbb{X}_r^o$ and $(p, A) \in \overline{J}^{2,+} u_\alpha(z)$.

By convexity of $\mathcal{H}$, (6.6) and the subsolution property of $u$, it is seen that

$$\mathcal{H}(p) \leq \alpha \mathcal{H}(p^*) + (1-\alpha) \mathcal{H}(D\psi(z)) \leq -(1-\alpha),$$

for all $(p, A) \in \overline{J}^{2,+} u_\alpha(z)$. This proves the second inequality in (6.10).

Recall that $\xi$ is defined via (6.9). By (6.5) $u_\alpha < u$ on $\overline{\mathbb{X}}_r$, and since $u = v = \phi$ on $\partial_r$, it is impossible that $\xi \in \partial_r$. Let $\eta = \eta(\xi)$ be as in Lemma 6.1. For $\gamma \in (1, \infty)$ and $\varepsilon \in (0,1)$, set

$$\begin{aligned} \Psi(x,y) &= |\gamma(x-y) - \varepsilon \eta|^2 + \varepsilon |y - \xi|^2, \\ \Phi(x,y) &= u_\alpha(x) - v(y) - \Psi(x,y), \qquad (x,y) \in \overline{\mathbb{X}}_r \times \overline{\mathbb{X}}_r, \end{aligned}$$

and let

$$(\widetilde{x}_{\varepsilon,\gamma}, \widetilde{y}_{\varepsilon,\gamma}) \equiv (\widetilde{x}, \widetilde{y}) \in \mathop{\arg\max}_{(x,y) \in \overline{\mathbb{X}}_r \times \overline{\mathbb{X}}_r} \Phi(x,y).$$

By Lemma 6.1,

(6.13) $$\xi + \frac{\varepsilon}{\gamma} \eta \in \mathbb{X}_r^o.$$

Clearly, $\Phi(\widetilde{x}, \widetilde{y}) \geq \Phi(\xi + \gamma^{-1}\varepsilon\eta, \xi)$. This can be rewritten as

$$(6.14) \quad u_\alpha(\widetilde{x}) - v(\widetilde{y}) - u_\alpha\left(\xi + \frac{\varepsilon}{\gamma}\eta\right) + v(\xi) \geq |\gamma(\widetilde{x} - \widetilde{y}) - \varepsilon\eta|^2 + \varepsilon|\widetilde{y} - \xi|^2.$$

Dividing by $\gamma^2$, we see that, for every $\varepsilon$, $|\widetilde{x} - \widetilde{y}| \to 0$ as $\gamma \to \infty$. This observation, along with (6.9), (6.14) and the continuity of $u_\alpha$ and $v$, gives that $\limsup_{\gamma \to \infty} |\gamma(\widetilde{x} - \widetilde{y}) - \varepsilon\eta|^2 + \varepsilon|\widetilde{y} - \xi|^2 \leq 0$. Hence, for all $\varepsilon \in (0, 1)$

$$(6.15) \qquad \widetilde{y} \to \xi, \qquad \gamma(\widetilde{x} - \widetilde{y}) \to \varepsilon\eta \qquad \text{as } \gamma \to \infty.$$

In particular, $\widetilde{x} = \widetilde{y} + \gamma^{-1}\varepsilon\eta + \gamma^{-1}o(1)$, as $\gamma \to \infty$. Hence, by (6.13) and Lemma 6.1, $\widetilde{x} \in \mathbb{X}_r^o$ for $\gamma > \gamma_0$, for some $\gamma_0 = \gamma_0(\varepsilon) < \infty$. By (6.5), (6.7), (6.9) and nonnegativity of $u$, it follows that $v(\xi) < u(\xi)$. By choosing $\gamma_0$ larger if necessary, we have $v(\widetilde{y}) < u(\widetilde{y})$ for $\gamma > \gamma_0$. Henceforth, assume $\gamma > \gamma_0$. Since $\xi \notin \partial_r$, we have that, for large $\gamma$,

$$(6.16) \qquad \widetilde{y} \in \mathbb{X}_r.$$

For $(x, r, p, A) \in \mathbb{X}_r \times \mathbb{R} \times \mathbb{R}^k \times \mathcal{S}(k)$, let

$$G(x, r, p, A) \doteq F(x, r, p, A) \vee \mathcal{H}(p).$$

Since $\widetilde{x} \in \mathbb{X}_r^o$, we have from (6.10) that

$$G(\widetilde{x}, u_\alpha(\widetilde{x}), p, X) \leq -(1 - \alpha), \qquad (p, X) \in \overline{J}^{2,+} u_\alpha(\widetilde{x}).$$

By (6.16) and the supersolution property of $v$,

$$G(\widetilde{y}, v(\widetilde{y}), q, Y) \geq 0 \qquad \forall (q, Y) \in \overline{J}^{2,-} v(\widetilde{y}).$$

Combining the above two displays and using the inequality $(a \vee b) - (c \vee d) \leq (a - c) \vee (b - d)$, $(a, b, c, d) \in \mathbb{R}^4$, we obtain that, for all $(p, X) \in \overline{J}^{2,+} u_\alpha(\widetilde{x})$ and $(q, Y) \in \overline{J}^{2,-} v(\widetilde{y})$,

$$1 - \alpha \leq G(\widetilde{y}, v(\widetilde{y}), q, Y) - G(\widetilde{x}, u_\alpha(\widetilde{x}), p, X)$$
$$(6.17) \qquad \leq \max\{|q - p|, \beta(v(\widetilde{y}) - u_\alpha(\widetilde{x})) + |b||q - p|$$
$$+ |\ell(\widetilde{x}) - \ell(\widetilde{y})| + (1/2)\text{trace}(\Sigma(X - Y))\}.$$

Next, noting that $u_\alpha(\widetilde{x}) - v(\widetilde{y}) \geq \Phi(\widetilde{x}, \widetilde{y}) \geq \Phi(\xi, \xi)$, and using (6.9), we have

$$(6.18) \qquad v(\widetilde{y}) - u_\alpha(\widetilde{x}) \leq \varepsilon^2 |\eta|^2.$$

Hence, by (6.17),

$$(6.19) \quad 1 - \alpha \leq \max\{|q - p|, \beta\varepsilon^2|\eta|^2 + |b||q - p|$$
$$+ |\ell(\widetilde{x}) - \ell(\widetilde{y})| + (1/2)\text{trace}(\Sigma(X - Y))\}.$$

By Theorem 3.2 of [8], for each $\sigma \in (0,\infty)$, one can find $X, Y \in \mathcal{S}(k)$ such that

$$(D_x\Psi(\widetilde{x},\widetilde{y}), X) \in \overline{J}^{2,+}u_\alpha(\widetilde{x}), \qquad (-D_y\Psi(\widetilde{x},\widetilde{y}), Y) \in \overline{J}^{2,-}v(\widetilde{y})$$

and

$$(6.20) \qquad \begin{pmatrix} X & 0 \\ 0 & -Y \end{pmatrix} \leq D^2\Psi(\widetilde{x},\widetilde{y}) + \sigma(D^2\Psi(\widetilde{x},\widetilde{y}))^2.$$

Observe that

$$(6.21) \qquad \begin{aligned} D_x\Psi(\widetilde{x},\widetilde{y}) &= 2\gamma(\gamma(\widetilde{x}-\widetilde{y})-\varepsilon\eta); \\ -D_y\Psi(\widetilde{x},\widetilde{y}) &= 2\gamma(\gamma(\widetilde{x}-\widetilde{y})-\varepsilon\eta) - 2\varepsilon(\widetilde{y}-\xi), \end{aligned}$$

$$(6.22) \qquad D^2\Psi(\widetilde{x},\widetilde{y}) = 2\gamma^2 \begin{pmatrix} I & -I \\ -I & I \end{pmatrix} + 2\varepsilon \begin{pmatrix} 0 & 0 \\ 0 & I \end{pmatrix}.$$

From (6.20), it follows that, with $\mathbf{1} = (I, I)'$,

$$X - Y = \mathbf{1}' \begin{pmatrix} X & 0 \\ 0 & -Y \end{pmatrix} \mathbf{1} \leq \mathbf{1}'(D^2\Psi(\widetilde{x},\widetilde{y}) + \sigma(D^2\Psi(\widetilde{x},\widetilde{y}))^2)\mathbf{1}.$$

This implies that

$$\text{trace}(\Sigma(X-Y)) \leq \text{trace}(\Sigma\mathbf{1}'(D^2\Psi(\widetilde{x},\widetilde{y}) + \sigma(D^2\Psi(\widetilde{x},\widetilde{y}))^2)\mathbf{1}).$$

From (6.22), it now follows that, for some $\varsigma \in (0,\infty)$ which is independent of $\varepsilon, \gamma, \sigma$,

$$(6.23) \qquad (1/2)\text{trace}(\Sigma(X-Y)) \leq \varsigma\varepsilon.$$

Using (6.18), (6.21) and (6.23) in (6.20), we have

$$1 - \alpha \leq \beta\varepsilon^2|\eta|^2 + 2\varepsilon|\widetilde{y}-\xi|(|b|+1) + |\ell(\widetilde{x}) - \ell(\widetilde{y})| + \varsigma\varepsilon.$$

Letting $\gamma \to \infty$, using (6.15) and continuity of $\ell$, we obtain

$$1 - \alpha \leq \beta\varepsilon^2|\eta|^2 + \varsigma\varepsilon.$$

Finally, letting $\varepsilon \to 0$ we arrive at a contradiction. Hence, (6.8) must hold, and the result follows. $\square$

**7. Uniqueness on unbounded domain.** This section proves uniqueness of solutions to (2.14) on $\mathbb{X}$. The two results stated below establish part (iii) of Theorem 2.1.

THEOREM 7.1. *Suppose that (2.10) holds. Let $u \in C_{\text{pol},+}(\mathbb{X})$ be a constrained viscosity solution of (2.14) on $\mathbb{X}$. Then $u \leq \overline{V}$ (in particular, $u \leq V$).*

THEOREM 7.2. *Suppose that* (2.10) *holds. Let* $u \in C_{\mathrm{pol},+}(\mathbb{X})$ *be a constrained viscosity solution of* (2.14) *on* $\mathbb{X}$. *Then* $u \geq V$ *provided that, in addition, one of the following holds:* (a) $u \in C_{\mathrm{pol},+}^{\mathrm{c}}(\mathbb{X})$; *or* (b) (2.11) *holds.*

While both results above hold regardless of the value of $m_\ell \in [0,\infty)$, in the special case where $m_\ell = 0$, the statement of part (a) of Theorem 7.2 is void: If a solution $u$ belongs to $C_{\mathrm{pol},+}$, then by Theorem 7.1, $u \leq V$, and by Lemma 4.4, $V$ is bounded, and so $u$ cannot lie in $C_{\mathrm{pol},+}^{\mathrm{c}}$.

Before proving the above, we show that these results, along with the results of the previous sections, imply Theorem 2.1.

PROOF OF THEOREM 2.1. For part (i), see Remark 4.2. Theorem 5.1 establishes part (ii). We now consider part (iii). Let $m_\ell > 0$ and $u \in C_{\mathrm{pol},+}^{\mathrm{c}}$ be a constrained viscosity solution of (2.14) on $\mathbb{X}$. Then Theorem 7.1 and Theorem 7.2(a) establish that $V = u$ and so by part (i) of the theorem, we get that $V$ is the only solution in the class $C_{\mathrm{pol},+}^{\mathrm{c}}$. For the case $m_\ell = 0$, Theorem 7.1 establishes the maximality of $V$ among solutions in $C_{\mathrm{pol},+}$. Finally, under (2.11), uniqueness in $C_{\mathrm{pol},+}$ follows from Theorem 7.1 and Theorem 7.2(b). This proves part (iii).

The identity $\overline{V} = V$ under condition (2.10) is an immediate consequence of parts (i)–(iii). □

PROOF OF THEOREM 7.1. Fix $u$ and let $c_u, m_u$ be such that $u(x) \leq c_u(1 + |x|^{m_u})$, $x \in \mathbb{X}$. Recall the notation of Lemmas 4.7 and 4.8. Fix $x \in \mathbb{X}$, a system $\Phi$, $Y \in \mathcal{A}_F(x)$ and $\varepsilon \in (0,1)$. Let $T^\varepsilon, Y^\varepsilon, X^\varepsilon, Z^\varepsilon$ be defined via (4.41), (4.42) and (4.43). Let $\sigma = \sigma(r)$ be as in (4.2). Now we estimate $\mathbb{E}[e^{-\beta\sigma}(|X_\sigma^\varepsilon|^{m_u} + 1)]$. Recalling the definition of $\bar{\mathcal{A}}_F(\Phi, x)$ and (4.12), we have that, for every $t$,

$$P(\sigma \leq t) = P(|\widehat{u}_1 \cdot X|_t^* \geq r) \leq 2a_0^{-1} r^{-1} \mathbb{E}(|X_t| + |B|_t^*) \to 0,$$

as $r \to \infty$. Hence, $\sigma(r) \to \infty$ a.s. as $r \to \infty$. By (4.12),

$$\mathbb{1}_{\sigma \leq T^\varepsilon} e^{-\beta\sigma}(|X_\sigma^\varepsilon|^{m_u} + 1) \leq c_0 \mathbb{1}_{\sigma \leq T^\varepsilon}[(|X_{T^\varepsilon}^\varepsilon| + |B|_{T^\varepsilon}^*)^{m_u} + 1].$$

Since $Y \in \bar{\mathcal{A}}_F(\Phi, x)$, $(|X_{T^\varepsilon}^\varepsilon| + |B|_{T^\varepsilon}^*)^{m_u}$ is integrable and, thus,

$$(7.1) \qquad \lim_{r \to \infty} \mathbb{E}[\mathbb{1}_{\sigma \leq T^\varepsilon} e^{-\beta\sigma}(|X_\sigma^\varepsilon|^{m_u} + 1)] = 0.$$

Let $g(\varepsilon) \equiv g_Y(\varepsilon) \doteq e^{-\beta T^\varepsilon} \mathbb{E}|X(T^\varepsilon)|^{m_u}$. Once more, since $Y \in \bar{\mathcal{A}}_F$, and $T^\varepsilon \to \infty$ as $\varepsilon \to 0$, we have $g(\varepsilon) \to 0$ as $\varepsilon \to 0$. Hence, by (4.41) and (4.44),

$$\mathbb{E}[\mathbb{1}_{\sigma > T^\varepsilon} e^{-\beta\sigma}(|X_\sigma^\varepsilon|^{m_u} + 1)]$$
$$\leq c_1(e^{-\beta T^\varepsilon}(\mathbb{E}|X_{T^\varepsilon}|^{m_u} + 1) + \mathbb{E}(\mathbb{1}_{\sigma > T^\varepsilon} e^{-\beta\sigma}|X_\sigma^\varepsilon - X_{T^\varepsilon}|^{m_u}))$$
$$\leq c_2(\varepsilon + g_Y(\varepsilon) + \mathbb{E}[\mathbb{1}_{\sigma > T^\varepsilon} e^{-\beta\sigma} \sup\{|B_s - B_{T^\varepsilon}|^{m_u} : s \in [T^\varepsilon, \sigma]\}])$$

$$= c_2\bigg(\varepsilon + g_Y(\varepsilon) + \sum_{i=0}^{\infty} \mathbb{E}[\mathbb{1}_{\sigma\in[T^\varepsilon+i,T^\varepsilon+i+1)}e^{-\beta(T^\varepsilon+i)}$$
$$\sup\{|B_s - B_{T^\varepsilon}|^{m_u} : s \in [T^\varepsilon, T^\varepsilon + i + 1]\}]\bigg)$$
$$\leq c_3\bigg(\varepsilon + g_Y(\varepsilon) + e^{-\beta T^\varepsilon}\sum_{i=0}^{\infty} e^{-\beta i}(1+i)^{m_u}\bigg)$$
$$\leq c_4(\varepsilon + g_Y(\varepsilon)).$$

Combining this with (7.1) we have

(7.2) $$\limsup_{r\to\infty} \mathbb{E}[e^{-\beta\sigma}(|X^\varepsilon_\sigma|^{m_u} + 1)] \leq c_5(\varepsilon + g_Y(\varepsilon)).$$

Next note that the control $Y^\varepsilon$ need not be in $\bar{\mathcal{A}}_r$, but we can modify it as follows. On $\{\sigma\infty\}$ let $X^{\varepsilon,r} = X^\varepsilon$. On $\{\sigma < \infty\}$, let $(X^{\varepsilon,r}(t), Y^{\varepsilon,r}(t)) = (X^\varepsilon(t), Y^\varepsilon(t))$ for $t \in [0, \sigma)$ and, recalling (4.17), set

$$X^{\varepsilon,r}(\sigma) = X^\varepsilon(\sigma-) + \gamma\Delta X^\varepsilon(\sigma), \qquad Y^{\varepsilon,r}(\sigma) = Y^\varepsilon(\sigma-) + \gamma\Delta Y^\varepsilon(\sigma),$$

where

$$\gamma = \gamma_r(X^\varepsilon(\sigma-), \Delta X^\varepsilon(\sigma)).$$

Then $X^{\varepsilon,r} = x + B + GY^{\varepsilon,r}$ on $[0,\sigma]$ and moreover, $Y^{\varepsilon,r} \in \bar{\mathcal{A}}_r$. Since $\gamma \leq 1$, $h(\Delta Y^{\varepsilon,r}(\sigma)) \leq h(\Delta Y^\varepsilon(\sigma))$, and by (2.2), $|X^{\varepsilon,r}(\sigma)| \leq a_0^{-1}\widehat{u}_1 \cdot X^{\varepsilon,r}(\sigma) \leq a_0^{-1}|X^\varepsilon(\sigma)|$. It therefore follows from Proposition 6.1 that

$$u(x) \leq \mathbb{E}\bigg[\int_{[0,\sigma]} e^{-\beta s}(\ell(X^{\varepsilon,r}_s)\,ds + h(dY^{\varepsilon,r}_s)) + c_u e^{-\beta\sigma}(|X^{\varepsilon,r}_\sigma|^{m_u} + 1)\bigg]$$
$$\leq \mathbb{E}\bigg[\int_{[0,\sigma]} e^{-\beta s}(\ell(X^\varepsilon_s)\,ds + h(dY^\varepsilon_s)) + c_u a_0^{-m_u} e^{-\beta\sigma}(|X^\varepsilon_\sigma|^{m_u} + 1)\bigg].$$

Thus, by (7.2) and Lemma 4.8,

$$u(x) \leq J(\Phi, x, Y^\varepsilon) + c_6(\varepsilon + g_Y(\varepsilon)) \leq J(\Phi, x, Y) + c_7(\varepsilon + g_Y(\varepsilon)).$$

Sending $\varepsilon \to 0$ and recalling that $Y \in \bar{\mathcal{A}}_F(\Phi, x)$ and $x \in \mathbb{X}$ are arbitrary, we conclude by Lemma 4.7 that $u \leq \overline{V}$ on $\mathbb{X}$. $\square$

PROOF OF THEOREM 7.2. Once again, let $u$ be fixed and let $c_u, m_u$ be such that $u(x) \leq c_u(1 + |x|^{m_u})$, $x \in \mathbb{X}$. Fix $x \in \mathbb{X}$. Let $\varepsilon \in (0, \varepsilon_0)$ be given, where the constant $\varepsilon_0 > 0$ will be chosen later. We will use the remarks below (2.13) and Proposition 6.1. We will construct a system $\widehat{\Phi} = (\widehat{\Omega}, \widehat{\mathcal{F}}, \widehat{\mathbb{P}}, \widehat{\mathcal{F}}_t, \widehat{B})$ and $\widehat{Y} \in \mathcal{A}(\widehat{\Phi}, x)$ such that $u(x) \geq J(x, \widehat{Y}) - c\varepsilon$ with a constant $c$ not depending on $\varepsilon$. This will clearly yield the result.

Consider first the case where $u \in C^{\mathrm{c}}_{\mathrm{pol},+}(\mathbb{X})$, namely, $u$ has compact level sets. Then $R_u(r) = \min_{\partial_r} u$ satisfies $R_u(r) \to \infty$ as $r \to \infty$. Let $\rho \in (0, \infty)$ and $\xi \in \mathbb{X}_\rho$. Consider the minimization problem associated with the right-hand side of (6.1) on $\mathbb{X}_r$, where $r > \rho$. If $Y \in \mathcal{A}_r(\xi)$ satisfies $J_{r,u}(\xi, Y) \leq V_{r,u}(\xi) + 1$, then from Proposition 6.1, with $\sigma = \sigma(r)$ [cf. (4.2)],

$$(7.3) \quad \mathbb{E}[e^{-\beta\sigma} u(X_\sigma)] + \mathbb{E} \int_{[0,\sigma]} e^{-\beta s} h(dY_s) \leq c_u(1 + |\xi|^{m_u}) + 1,$$

where $X$ is the controlled process corresponding to $Y$. Since $X_\sigma \in \partial_r$, we have

$$(7.4) \quad u(X_{\sigma(r)}) \geq R_u(r).$$

Therefore, by (7.3),

$$\mathbb{P}(\sigma(r) < 1) \leq e^\beta R_u(r)^{-1} [c_u(1 + |\xi|^{m_u}) + 1].$$

This shows the following:

$$(7.5) \quad \begin{array}{c} \text{For every } \rho, \text{ one can find } r = \overline{r}(\rho) > \rho \\ \text{such that } \mathbb{P}(\sigma(r) < 1) \leq 1/2, \ \forall \xi \in \mathbb{X}_\rho. \end{array}$$

Next consider the case where (2.11) is assumed and $u \in C_{\mathrm{pol},+}$. By (2.6) and (2.11),

$$\int_{[0,\sigma]} e^{-\beta s} h(dY_s) \geq c_h \int_{[0,\sigma]} e^{-\beta s} d|Y|_s \geq c_h e^{-\beta \sigma} |Y_\sigma|.$$

Since $\widehat{u}_1 \cdot \xi \leq \rho$ and $\widehat{u}_1 \cdot X_\sigma r$, we have that

$$\mathbb{E} \int_{[0,\sigma]} e^{-\beta s} h(dY_s) \geq c_h |G|^{-1} \mathbb{E} e^{-\beta\sigma} \widehat{u}_1 \cdot GY_\sigma$$

$$(7.6) \quad \geq c_h |G|^{-1} \mathbb{E}[e^{-\beta\sigma}(r - \rho - \widehat{u}_1 \cdot B_\sigma)]$$

$$\geq c_h |G|^{-1} (r - \rho) \mathbb{E} e^{-\beta\sigma} - c_h |G|^{-1} \mathbb{E}[e^{-\beta\sigma} |B|_\sigma^*],$$

where $|G| > 0$ denotes the operator norm of $G$. Note that

$$\mathbb{E}[e^{-\beta\sigma} |B|_\sigma^*] \leq \sum_{i=0}^\infty \mathbb{1}_{\sigma \in [i, i+1)} e^{-\beta i} c_1(i+1) \leq c_2,$$

where $c_2$ does not depend on $Y$ and $r$. Thus,

$$\mathbb{E} \int_{[0,\sigma]} e^{-\beta s} h(dY_s) \geq c_h |G|^{-1} [(r-\rho) \mathbb{E} e^{-\beta\sigma} - c_2].$$

Combining this with (7.3), we have that (7.5) holds in this case as well.

Define inductively a sequence of domains $\mathbb{X}_{r_n}, \mathbb{X}_{\rho_n}, \widehat{u}_1 \cdot x = \rho_0 < r_1 < \rho_1 < r_2 < \cdots$ via the relations $r_n = \overline{r}(\rho_{n-1})$, $n \in \mathbb{N}$, and $\rho_n = r_n + 1$, $n \in \mathbb{N}$. For $n \in \mathbb{N}$ and $\delta > 0$, denote

$$m_n(\delta) = \max\{u(\xi) - u(z) : \xi, z \in \mathbb{X}_{\rho_n}, |\xi - z| \leq \delta\}.$$

By assumption, $\widehat{u}_0 \in \mathbb{U}^o$ [cf. the comment following (2.1)]. Hence, there is a constant $a_1 > 0$ such that $B_{a_1}(\widehat{u}_0) \subset \mathbb{U}^o \cap \mathbb{X}^o$. Fix such $a_1$. Fix also $k$ linearly independent unit vectors $u_i \in B_{a_1}(\widehat{u}_0)$ such that, with $\mathbb{A} \doteq \text{cone}(u_i, i = 1, \ldots, k)$, one has $\widehat{u}_0 \in \mathbb{A}^o$. Thus, there exists $a > 0$ such that $B_a(\widehat{u}_0) \subset \mathbb{A}$. Let such $a$ be fixed and denote $c_3 = (1 + a^{-1})k^{1/2}$. Let $\varrho : \mathbb{A} \to \mathbb{Y}$ be the linear map of Lemma A.2 and denote by $|\varrho|$ its operator norm. For $n \in \mathbb{N}$, let $\mathfrak{S}_n$ denote the finite set

$$\mathfrak{S}_n = (\lambda_{n,\varepsilon}\mathbb{Z}^k) \cap \mathbb{X}_{\rho_n},$$

where $\lambda_{n,\varepsilon} > 0$ are fixed constants, so small that, for $n \in \mathbb{N}$,

(7.7) $\quad (1 + |\varrho|c_3)\lambda_{n,\varepsilon} \leq \varepsilon 2^{-n} \quad \text{and} \quad m_n(c_3\lambda_{n,\varepsilon}) \leq \varepsilon 2^{-n}.$

We next show that

(7.8) $\quad \forall \xi \in \partial_{r_n} \ \exists z \in \mathfrak{S}_n \ \text{s.t.} \ z - \xi \in \mathbb{A}, \quad |z - \xi| \leq c_3\lambda_{n,\varepsilon}.$

Indeed, given $\xi \in \partial_{r_n}$ we have $\xi + C \subset \mathbb{X}$, where $C \subset \mathbb{A}$ denotes the cone generated by $B_a(\widehat{u}_0)$. Since $C$ can be written as $\bigcup_{\alpha > 0} B_{\alpha a}(\alpha \widehat{u}_0)$, we have

$$\mathcal{Z} \doteq \xi + B_{\lambda_{n,\varepsilon}k^{1/2}}(a^{-1}\lambda_{n,\varepsilon}k^{1/2}\widehat{u}_0) \subset \mathbb{X}.$$

It is easy to see that, for every $w \in \mathbb{R}^k$, there exists $z \in \lambda_{n,\varepsilon}\mathbb{Z}^k$ such that $|z - w| < k^{1/2}\lambda_{n,\varepsilon}$. This shows that $\mathcal{Z}$ contains a point in $\lambda_{n,\varepsilon}\mathbb{Z}^k$. Choosing $\varepsilon_0 = c_3^{-1}$, we have $|z - \xi| < c_3\lambda_{n,\varepsilon} < 1$ for every $z \in \mathcal{Z}$ and $\varepsilon \in (0, \varepsilon_0)$, and thus for all such $\varepsilon$, $\mathcal{Z} \subset \mathbb{X}_{r_n+1}$. Thus, $\mathcal{Z}$ contains a point $z \in \mathfrak{S}_n$. By construction, $z - \xi \in \mathbb{A}$ and, therefore, (7.8) holds.

For $n \in \mathbb{N}$, let $M_n : \partial_{r_n} \to \mathfrak{S}_n$ denote a measurable map such that, for every $\xi \in \partial_{r_n}$, condition (7.8) is met by $z = M_n(\xi)$. By Lemma A.2,

(7.9) $\quad G\varrho(M_n(\xi) - \xi) = M_n(\xi) - \xi.$

Consider a complete probability space $(\widehat{\Omega}, \widehat{\mathcal{F}}, \widehat{\mathbb{P}})$ supporting countably many independent $(b, \Sigma)$ Brownian motions. In particular, let a $(b, \Sigma)$ Brownian motion $B^{n,z}$ be associated with each $n \in \mathbb{N}$ and $z \in \mathfrak{S}_n$, and let a $(b, \Sigma)$ Brownian motion $B^{0,x}$ be associated with $x$. For each $z \in \mathfrak{S}_n$, consider the minimization problem associated with the right-hand side of (6.1), substituting $r_{n+1}$ for $r$ and $z$ for $x$. We write the system $(\widehat{\Omega}, \widehat{\mathcal{F}}, \widehat{\mathcal{F}}_t^{n,z}, \widehat{\mathbb{P}}, B^{n,z})$, where $\widehat{\mathcal{F}}_t^{n,z}$ is the $\widehat{\mathbb{P}}$ completion of the $\sigma$-field generated by $B^{n,z}$, as $\widehat{\Phi}^{n,z}$. $\widehat{\Phi}^{0,x}$ is defined similarly.

Using Proposition 6.1, find a $Y^{n,z} \in \mathcal{A}_{r_{n+1}}(z, \widehat{\Phi}^{n,z})$, for which

(7.10) $\quad \mathbb{E}\left[\int_{[0,\sigma(r_{n+1})]} e^{-\beta t}(\ell(X_t^{n,z})\,dt + h(dY_t^{n,z})) + e^{-\beta\sigma(r_{n+1})}u(X_{\sigma(r_{n+1})}^{n,z})\right]$
$\qquad \leq u(z) + \varepsilon 2^{-n},$

where $X^{n,z}$ is the controlled process corresponding to $Y^{n,z}$. To account for the dependence of $\sigma(r_{n+1})$ on the initial point $z \in \mathfrak{S}_n$, it will be more convenient to write in what follows $\sigma^{n,z}$ for $\sigma(r_{n+1})$. By construction and by (7.5), for each pair $(n,z)$,

$$\mathbb{P}(\sigma^{n,z} < 1) \leq 1/2. \tag{7.11}$$

Define inductively a sequence of processes $(\widehat{B}_n, \widehat{Y}_n, \widehat{X}_n)$ as follows. Let $\widehat{B}_1 = B^{0,x}$, and $\widehat{\sigma}_1 = \sigma^{0,x}$. Let $\Xi_1 = X^{0,x}(\sigma^{0,x})$ and note that $\Xi_1 \in \partial_{r_1}$. Let also

$$Z_1 \doteq M_1(\Xi_1), \qquad \widehat{Y}_1 \doteq Y^{0,x} \mathbb{1}_{[0,\sigma^{0,x}]} + \varrho(Z_1 - \Xi_1) \mathbb{1}_{(\sigma^{0,x},\infty)}, \tag{7.12}$$

and $\widehat{X}_1 = x + \widehat{B}_1 + G\widehat{Y}_1$ on $[0, \widehat{\sigma}_1]$ [we do not define $\widehat{X}_1$ on $(\widehat{\sigma}_1, \infty)$]. Note that by (7.9) and (7.12), $\widehat{X}_1(\widehat{\sigma}_1) = Z_1$, so by (7.8), $|\widehat{X}_1(\widehat{\sigma}_1) - \Xi_1| \leq c_3 \lambda_{1,\varepsilon}$. Also note that $\Delta \widehat{Y}_1(\widehat{\sigma}_1) = \Delta Y^{0,x}(\widehat{\sigma}_1) + \varrho(Z_1 - \Xi_1)$. Hence, by (4.6), (7.7) and (7.10)

$$\mathbb{E}\left[\int_{[0,\widehat{\sigma}_1]} e^{-\beta t}(\ell(\widehat{X}_1(t))\,dt + h(d\widehat{Y}_1(t))) + e^{-\beta\widehat{\sigma}_1} u(\widehat{X}_1(\widehat{\sigma}_1))\right]$$
$$\leq u(x) + \varepsilon 2^{-1} + |h||\varrho|c_3\lambda_{1,\varepsilon} + m_1(c_3\lambda_{1,\varepsilon}) \tag{7.13}$$
$$\leq u(x) + c_4 \varepsilon 2^{-1}.$$

Consider $n \geq 2$. On the set $\{\widehat{\sigma}_{n-1} = \infty\}$, let $\widehat{\sigma}_n = \infty$ and $(\widehat{B}_n, \widehat{Y}_n, \widehat{X}_n) = (\widehat{B}_{n-1}, \widehat{Y}_{n-1}, \widehat{X}_{n-1})$. Next consider the set $\{\widehat{\sigma}_{n-1} < \infty\}$. Let

$$\widehat{B}_n = \widehat{B}_{n-1} \mathbb{1}_{[0,\widehat{\sigma}_{n-1}]} + (B^{n-1,Z_{n-1}} - B^{n-1,Z_{n-1}}(\widehat{\sigma}_{n-1})) \mathbb{1}_{(\widehat{\sigma}_{n-1},\infty)}$$

and $\widehat{\sigma}_n = \widehat{\sigma}_{n-1} + \sigma^{n-1,Z_{n-1}}$. Let $\Xi_n = X^{n-1,Z_{n-1}}(\sigma^{n-1,Z_{n-1}})$, and

$$Z_n = M_n(\Xi_n), \qquad \widehat{Y}_n = \widehat{Y}_{n-1} \mathbb{1}_{[0,\widehat{\sigma}_n]} + \varrho(Z_n - \Xi_n) \mathbb{1}_{(\widehat{\sigma}_n,\infty)}.$$

Let $\widehat{X}_n = x + \widehat{B}_n + \widehat{Y}_n$ on $[0, \widehat{\sigma}_n]$ [and we have not defined it on $(\widehat{\sigma}_n, \infty)$]. By (7.8), $|\widehat{X}_n(\widehat{\sigma}_n) - \Xi_n| \leq c_3 \lambda_{n,\varepsilon}$. Denoting the filtration generated by $\widehat{B}_n$ as $\widehat{\mathcal{F}}_t^n$, we have from (7.10), in a manner similar to the proof of (7.13), that

$$\mathbb{E}\left[\int_{(\widehat{\sigma}_{n-1},\widehat{\sigma}_n]} e^{-\beta t}(\ell(\widehat{X}_n(t))\,dt + h(d\widehat{Y}_n(t))) + e^{-\beta\widehat{\sigma}_n} u(\widehat{X}_n(\widehat{\sigma}_n)) \Big| \widehat{\mathcal{F}}^n_{\widehat{\sigma}_{n-1}}\right]$$
$$\leq u(\widehat{X}_n(\widehat{\sigma}_{n-1})) + c_4 \varepsilon 2^{-n},$$

a.s. Iterating the above inequality and using (7.13), we now have

$$\mathbb{E}\left[\int_{[0,\widehat{\sigma}_n]} e^{-\beta t}(\ell(\widehat{X}_n(t))\,dt + h(d\widehat{Y}_n(t))) + e^{-\beta\widehat{\sigma}_n} u(X_n(\widehat{\sigma}_n))\right]$$
$$\leq u(x) + c_4 \sum_{i=1}^n \varepsilon 2^{-n} \tag{7.14}$$
$$\leq u(x) + c_4 \varepsilon.$$

To see that $\widehat{\sigma}_n \to \infty$ a.s., let $F_0 = \{\varnothing, \Omega\}$, and for $n \in \mathbb{N}$, let $F_n$ be the sigma-field generated by $F_{n-1}$ and $(\widehat{B}_n(s), \widehat{Y}_n(s) : s \in [0, \infty))$. Then $\widehat{\sigma}_n - \widehat{\sigma}_{n-1} \in F_n$ and, by (7.11), $\mathbb{P}(\widehat{\sigma}_n - \widehat{\sigma}_{n-1} \geq 1 | F_{n-1}) > 1/2$. By the second Borel–Cantelli lemma (cf. [13], page 240), $\widehat{\sigma}_n \to \infty$ a.s.

Since $\widehat{\sigma}_n \to \infty$ a.s., the limits $\widehat{B} = \lim_n \widehat{B}_n$, $\widehat{Y} = \lim_n \widehat{Y}_n$ are well defined outside a null set. Let $\widehat{X} = x + G\widehat{Y} + \widehat{B}$. By construction, the process $\widehat{B}$ is a $(b, \Sigma)$-Brownian motion, and $\widehat{Y} \in \mathcal{A}(x, \widehat{\Phi}, \widehat{B})$. Finally, by (7.14), $J(x, \widehat{Y}) \leq u(x) + c_4 \varepsilon$. Hence, $V(x) \leq u(x) + c_4 \varepsilon$, and the result follows since $\varepsilon > 0$ is arbitrarily. $\square$

**8. Dynamic programming principles.** In this section we prove Propositions 5.1, 5.2 and Lemma 5.1.

PROOF OF PROPOSITION 5.1. Fix $t \in (0, \infty)$, $x \in \mathbb{X}$ and $Y \in \mathcal{A}(x)$. For brevity, we denote $\tau \wedge t$ by $\theta$. Note that

$$(8.1) \quad \mathbb{E}\left[\int_{[\theta,\infty)} e^{-\beta s}(\ell(X_s)\,ds + h(dY_s^1))|\mathcal{F}_\theta\right] \geq e^{-\beta\theta} V(X_\theta),$$

where $Y_s^1 \doteq Y_s \mathbb{1}_{s<\theta} + (Y_s - \Delta Y_\theta)\mathbb{1}_{s \geq \theta}$. The proof of (8.1) follows in a straightforward manner on recalling that $\overline{\mathcal{F}}_t$ is generated by $W$ and using the strong Markov property of $W$. This immediately shows that

$$J(x,Y) \geq \mathbb{E}\left[\int_{[0,\theta)} e^{-\beta s}(\ell(X_s)\,ds + h(dY_s)) + e^{-\beta\theta}(h(\Delta Y_\theta) + V(X_\theta))\right].$$

Taking infimum over $Y \in \mathcal{A}(x)$ in the above inequality, we have that

$$(8.2) \quad V(x) \geq \inf_{Y \in \mathcal{A}(x)} \mathbb{E}\left[\int_{[0,\theta]} e^{-\beta s}(\ell(X_s)\,ds + h(dY_s)) + e^{-\beta\theta}(V(X_\theta))\right].$$

Now we prove the reverse inequality. Once again, fix $x \in \mathbb{X}$ and $Y \in \mathcal{A}(x)$. Let $\delta \in (0, \infty)$ be arbitrary. Then

$$(8.3) \quad \begin{array}{c} \exists \widetilde{Y} \in \mathcal{A}(x) \text{ and the corresponding } \widetilde{X} \text{ s.t.} \\ \widetilde{Y}(s) = Y(s), s \in [0, \theta) \text{ and } (8.4) \text{ holds:} \end{array}$$

$$(8.4) \quad e^{-\beta\theta} V(X_\theta) \geq \mathbb{E}\left[\int_{[\theta,\infty)} e^{-\beta s}(\ell(\widetilde{X}_s)\,ds + h(d\widetilde{Y}_s)) - e^{-\beta\theta} h(\Delta Y_\theta)|\mathcal{F}_\theta\right] - \delta.$$

The proof of (8.3) is provided in the Appendix. This shows that

$$\mathbb{E}\left[\int_{[0,\theta]} e^{-\beta s}(\ell(X_s)\,ds + h(dY_s)) + e^{-\beta\theta}(V(X_\theta))\right]$$
$$\geq \mathbb{E}\int_{[0,\infty)} e^{-\beta s}(\ell(\widetilde{X}_s)\,ds + h(d\widetilde{Y}_s)) - \delta$$
$$\geq V(x) - \delta.$$

Taking infimum over all $Y \in \mathcal{A}(x)$, the above inequality and (8.2) establish the result. □

PROOF OF PROPOSITION 5.2. The proof is adapted from that of Proposition III.1.1 of [4]. In view of Lemma 5.1, it suffices to work with $\widehat{V}$. Let $\widetilde{\Phi}$ be a system, $\zeta$ be a $\widetilde{\mathcal{F}}_0$ measurable $\mathbb{R}^k$ valued random variable with probability distribution $\mu$, $\widetilde{Y} \in \widehat{\mathcal{A}}(\widetilde{\Phi}, \zeta)$ and $\widetilde{X}_t = \zeta + B_t + G\widetilde{Y}_t$ be the corresponding controlled process. From arguments in Theorem I.1.6 of [4], it follows that the conditional law of $\widetilde{X}$ given $\zeta = x$ is $\mu$-almost surely the law of a controlled process $\widetilde{X}^x$ corresponding to some $\widetilde{Y}^x \in \widehat{\mathcal{A}}(\widetilde{\Phi}^x, x)$ and some system $\widetilde{\Phi}^x$. Furthermore,

$$\int J(x, \widetilde{\Phi}^x, \widetilde{Y}^x) \mu(dx) = J(\widetilde{\Phi}, \zeta, Y).$$

Also, following Lemma III.1.1 of [4], we have that, given $\varepsilon > 0$, there exists an extended system $(\widetilde{\Phi}, \zeta)$ and $\widetilde{Y} \in \widehat{\mathcal{A}}(\Phi, \zeta)$ such that $\widetilde{Y}^x$ defined as above is $\varepsilon$-optimal for $\widehat{V}(x)$, for $\mu$-a.e. $x$. That is

$$J(x, \widetilde{\Phi}^x, \widetilde{Y}^x) \le \widehat{V}(x) + \varepsilon, \qquad \mu\text{- a.e. } x.$$

We call such a $\widetilde{Y}$ an $\varepsilon$-optimal control for $(\widetilde{\Phi}, \zeta)$. Now fix $x \in \mathbb{X}$. Let $Y^1 \in \widehat{\mathcal{A}}(\Phi, x)$ and let $X^1$ be the corresponding controlled process and $\tau = \tau_{Y^1}$ be defined via (5.1). Once more we will denote $t \wedge \tau$ by $\theta$. Let $\mu$ now denote the probability distribution of $X^1(\theta)$ and let $\widetilde{Y} \in \widehat{\mathcal{A}}(\widetilde{\Phi}, \zeta)$ be an $\varepsilon$-optimal control given on some extended system $(\widetilde{\Phi}, \zeta)$, with probability law of $\zeta$ equal to $\mu$. Let $\widetilde{X}$ be the corresponding controlled process. By augmenting $\Phi$ suitably, one can construct on it processes $X$ and $Y$ such that $Y \in \widehat{\mathcal{A}}(\Phi, x)$, $X = x + B + GY$, $X^1(s) = X(s)$, $Y^1(s) = Y(s)$ for $s \in [0, \theta)$ and the conditional distribution of $(X(\theta + \cdot), Y(\theta + \cdot) - Y^1(\theta))$ given $\widehat{\mathcal{F}}_\theta$ is, for almost every $\omega$, the same as the distribution of $(\widetilde{X}^z(\cdot), \widetilde{Y}^z(\cdot))$, with $z = X^1(\theta(\omega), \omega)$. Here $\widehat{\mathcal{F}}_t = \sigma\{X^1(s), Y^1(s), s \le t\}$. For details on this construction, we refer the reader to the proof of Theorem III.1.1 of [4]. Next, setting $\widehat{Y}(s) \doteq Y(s+\theta) - Y^1(\theta)$, we have

$$\widehat{V}(x) \le \mathbb{E}\left[\int_{[0,\theta]} e^{-\beta s}(\ell(X_s)\,ds + h(dY_s)) + \int_{(\theta,\infty)} e^{-\beta s}(\ell(X_s)\,ds + h(dY_s))\right]$$

$$\le \mathbb{E}\left[\int_{[0,\theta]} e^{-\beta s}(\ell(X_s^1)\,ds + h(dY_s^1)) + e^{-\beta\theta}h(\widehat{Y}(0))\right]$$

$$+ \mathbb{E}\left[e^{-\beta\theta}\mathbb{E}\left[\int_{(0,\infty)} e^{-\beta s}(\ell(X_{\theta+s})\,ds + h(d\widehat{Y}_s))|\widehat{\mathcal{F}}_\theta\right]\right]$$

$$= \mathbb{E}\int_{[0,\theta]} e^{-\beta s}(\ell(X_s^1)\,ds + h(dY_s^1)) + \mathbb{E}[e^{-\beta\theta}J(x, \widetilde{\Phi}^x, \widetilde{Y}^x)|_{x=X^1(\theta)}]$$

$$\leq \mathbb{E}\left[\int_{[0,\theta]} e^{-\beta s}(\ell(X_s^1)\,ds + h(dY_s^1)) + e^{-\beta\theta}\widehat{V}(X_\theta^1)\right] + \varepsilon,$$

where the second inequality follows on observing that $h(\Delta Y_\theta) \leq h(\Delta Y_\theta^1) + h(\widehat{Y}(0))$. Taking infimum over all $Y^1 \in \widehat{\mathcal{A}}(\Phi, x)$ and over all systems $\Phi$ and letting $\varepsilon \to 0$, we obtain

(8.5)
$$\widehat{V}(x) \leq \inf_\Phi \inf_{Y^1 \in \widehat{\mathcal{A}}(\Phi_1, x)} \mathbb{E}\left[\int_{[0,\theta]} e^{-\beta s}(\ell(X_s^1)\,ds + h(dY_s^1)) + e^{-\beta\theta}\widehat{V}(X_\theta^1)\right].$$

Next, let $\Phi^1$ be a system and $Y^1 \in \widehat{\mathcal{A}}(\Phi^1, x)$ be $\varepsilon$-optimal for $\widehat{V}(x)$. Then

(8.6)
$$\widehat{V}(x) + \varepsilon \geq \mathbb{E}\int_{[0,\theta]} e^{-\beta s}(\ell(X_s^1)\,ds + h(dY_s^1))$$
$$+ \mathbb{E}\int_{(\theta,\infty)} e^{-\beta s}(\ell(X_s^1)\,ds + h(dY_s^1)).$$

Once more, the arguments in Theorem I. 1.6 of [4] yield that the probability law of $(X^1(\theta+\cdot), Y^1(\theta+\cdot) - Y^1(\theta))$ given $\widehat{\mathcal{F}}_\theta$ is, for almost all $\omega$, the law of $(X_x^1, Y_x^1)$, where $Y_x^1 \in \widehat{\mathcal{A}}(\Phi^x, x)$, $\Delta Y_x^1(0) = 0$ a.s., $\Phi^x$ is some system, $X_x^1$ is the corresponding controlled process and $x = X^1(\theta(\omega), \omega)$. This shows that

$$\mathbb{E}\left[\int_{(\theta,\infty)} e^{-\beta s}(\ell(X_s^1)\,ds + h(dY_s^1))|\widehat{\mathcal{F}}_\theta\right] \geq e^{-\beta\theta} J(x, \Phi^x, Y_x^1)|_{x=X^1(\theta(\omega),\omega)}.$$

Thus,

$$\mathbb{E}\int_{(\theta,\infty)} e^{-\beta s}(\ell(X_s^1)\,ds + h(dY_s^1)) \geq \mathbb{E}[e^{-\beta\theta}\widehat{V}(X_\theta^1)].$$

Substituting the above in (8.7) and taking $\varepsilon \to 0$ yields

(8.7)
$$\widehat{V}(x) \geq \inf_\Phi \inf_{Y^1 \in \widehat{\mathcal{A}}(\Phi_1, x)} \mathbb{E}\left[\int_{[0,\theta]} e^{-\beta s}(\ell(X_s^1)\,ds + h(dY_s^1)) + e^{-\beta\theta}\widehat{V}(X_\theta^1)\right].$$

The result follows on combining (8.5) and (8.7). □

PROOF OF LEMMA 5.1. Fix $x \in \mathbb{X}$, a system $\Phi$ and $Y \in \bar{\mathcal{A}}(\Phi, x)$ such that $J(x, \Phi, Y) < \infty$. In view of Lemmas 4.7 and 4.8, we can assume without loss of generality [cf. (4.44)] that

(8.8) $\quad e^{-\beta t}\mathbb{E}|GY_t| \to 0 \quad \text{as } t \to \infty \text{ and } \mathbb{E}\int_{[0,\infty)} e^{-\beta t}|GY_t|\,dt < \infty.$

Write $Y_t = Y_t^c + Y_t^d$, where $Y_t^d \doteq \sum_{0 \leq s \leq t} \Delta Y_s$ and $Y_t^c = Y_t - Y_t^d$. Note that both $Y^c$ and $Y^d$ are RCLL, $\overline{\mathcal{F}}_t$ adapted and have increments in $\mathbb{Y}$. Furthermore, the measures $d|Y^c|$ and $d|Y^d|$ are singular and, therefore, from Lemma 4.2,

$$J(\Phi, x, Y) = \mathbb{E} \int_{[0,\infty)} e^{-\beta s}[\ell(X_s)\, ds + h(dY_s^d) + h(dY_s^c)],$$

where $X$ is the controlled process corresponding to $Y$. Fix $\delta > 0$ and define

$$Y_t^\delta \doteq Y_{(k-2)\delta}^c + \frac{1}{\delta}(t - (k-1)\delta)(Y_{(k-1)\delta}^c - Y_{(k-2)\delta}^c),$$
(8.9)
$$t \in [(k-1)\delta, k\delta), k = 2, 3, \ldots.$$

Set $Y_t^\delta = Y_0^c = 0$ for $t \in [0, \delta)$. Note that $Y^\delta$ is continuous $\overline{\mathcal{F}}_t$ adapted and has increments in $\mathbb{Y}$. Also let

(8.10) $X_t^\delta \doteq \Gamma(x + B + GY^\delta + GY^d)(t) = x + B_t + GY_t^\delta + GY_t^d + G\widehat{y}_0 \vartheta_t^\delta,$

where $\vartheta_t^\delta \doteq \widehat{\Gamma}(x + B + GY^\delta + GY^d)(t)$. Define $\xi_t^\delta \doteq X_t^\delta - G\widehat{y}_0 \vartheta_t^\delta = x + B_t + GY_t^d + GY_t^\delta$. Note that $\Delta \xi_t^\delta = G\Delta Y_t^d = G\Delta Y_t$. Thus, $\Delta \xi_t^\delta = G\mathbb{E}(\Delta Y_t | \Delta \xi_t^\delta)$. Define $\widetilde{Y}_t^d \doteq \sum_{0 \leq s \leq t} \mathbb{E}(\Delta Y_s | \Delta \xi_s^\delta)$. Clearly, $\widetilde{Y}^d$ is adapted to the filtration $\sigma\{X_s^\delta, \vartheta_s^\delta, s \leq t\}$. Furthermore, using the convexity of $h$, it is easy to see that

(8.11) $$\mathbb{E} \int_{[0,\infty)} e^{-\beta t} h(d\widetilde{Y}_t^d) \leq \mathbb{E} \int_{[0,\infty)} e^{-\beta t} h(dY_t^d).$$

Also, it is easy to check that $G\widetilde{Y}_t^d = GY_t^d$. Now let $\eta_t^\delta \doteq x + B_t + GY_t^\delta = x + B_t + \int_0^t G\dot{Y}_s^\delta\, ds$, where $\dot{Y}_s^\delta \doteq \frac{1}{\delta}(Y_{(k-1)\delta}^c - Y_{(k-2)\delta}^c)\mathbb{1}_{s \in [(k-1)\delta, k\delta)}$ and $\dot{Y}_s^\delta = 0$ for $s \in [0, \delta)$. Note that $\eta_t^\delta = \xi_t^\delta - G\widetilde{Y}_t^\delta$ and, therefore, it is $\widehat{\mathcal{F}}_t$ adapted, where $\widehat{\mathcal{F}}_t = \sigma\{X_s^\delta, \vartheta_s^\delta : s \leq t\}$. Thus, from Theorem 4.2 of [33], there exists a $(b, \Sigma) - \widehat{\mathcal{F}}_t$-Brownian motion $\widetilde{B}_t$ such that

(8.12) $$\eta_t^\delta = x + \widetilde{B}_t + G \int_0^t \mathbb{E}(\dot{Y}_s^\delta | \widehat{\mathcal{F}}_s)\, ds.$$

Note that $\widetilde{Y}_t^\delta \doteq \int_0^t \mathbb{E}(\dot{Y}_s^\delta | \widehat{\mathcal{F}}_s)\, ds$ is continuous, $\widehat{\mathcal{F}}_t$ adapted and has increments in $\mathbb{Y}$. Once more, using the radial homogeneity and convexity of $h$, we have

(8.13) $$\mathbb{E} \int_{[0,\infty)} e^{-\beta t} h(d\widetilde{Y}_t^\delta) \leq \mathbb{E} \int_{[0,\infty)} e^{-\beta t} h(dY_t^\delta).$$

Thus, defining $\widehat{Y}_t^\delta \doteq \widetilde{Y}_t^\delta + \widetilde{Y}_t^d$, we have by (8.10) and (8.12) $X_t^\delta = x + \widetilde{B}_t + G\widehat{Y}_t^\delta + G\widehat{y}_0 \vartheta_t^\delta$ and therefore, $\overline{Y}_t^\delta \doteq \widehat{Y}_t^\delta + \widehat{y}_0 \vartheta_t^\delta \in \widehat{\mathcal{A}}(\widetilde{\Phi}, x)$, where $\widetilde{\Phi} = (\Omega, \mathcal{F}, (\overline{\mathcal{F}}_t), \mathbb{P}, \widetilde{B})$. In order to prove the proposition, it now suffices to show that

(8.14) $$\limsup_{\delta \to 0} J(x, \widetilde{\Phi}, \overline{Y}^\delta) \leq J(x, \Phi, Y).$$

From (8.10) it follows that, for all $T \in (0, \infty)$,

$$\sup_{0 \leq t \leq T}[|X^\delta(t) - X(t)| + \vartheta^\delta(t)] \leq c_1 \sup_{0 \leq t \leq T}|GY^\delta(t) - GY^c(t)|$$

(8.15)

$$\leq c_1 \sup_{0 \leq t \leq T}|Y^\delta(t) - Y^c(t)|$$

and due to the sample path continuity of $Y^c$, the last term approaches 0 as $\delta \to 0$. This, along with continuity of $\ell$ and the Lipschitz continuity of $h$, shows that, for all $t \in (0, \infty)$,

$$\ell(X_t^\delta) \to \ell(X_t) \quad \text{and} \quad \sup_{0 \leq s \leq t}(\vartheta_s^\delta + |h(Y_s^\delta) - h(Y_s^c)|) \to 0$$

(8.16)

a.s. as $\delta \to 0$.

Also, using (2.8) and (2.2), it can be seen that

$$\ell(X_s^\delta) \leq c_2(1 + (|GY|_s^*)^{m_\ell} + (|GY^c|_s^*)^{m_\ell} + (|B|_s^*)^{m_\ell})$$
$$\leq c_3(1 + \ell(X_s) + (|B|_s^*)^{m_\ell}).$$

Thus, recalling that $\mathbb{E}\int_{[0,\infty)} e^{-\beta t}\ell(X_s)\,ds < J(x, \Phi, Y) < \infty$, we have from (8.16) and an application of dominated convergence theorem that, as $\delta \to 0$,

(8.17) $$\mathbb{E}\int_{[0,\infty)} e^{-\beta t}\ell(X_s^\delta)\,ds \to \mathbb{E}\int_{[0,\infty)} e^{-\beta t}\ell(X_s)\,ds.$$

Next, combining (8.11) and (8.13), we have

(8.18)
$$\mathbb{E}\int_{[0,\infty)} e^{-\beta s}h(d\overline{Y}_s^\delta) \leq \mathbb{E}\int_{[0,\infty)} e^{-\beta s}h(d\widetilde{Y}_s^d) + \mathbb{E}\int_{[0,\infty)} e^{-\beta s}h(d\widetilde{Y}_s^\delta)$$
$$+ h(G\widehat{y}_0)\mathbb{E}\int_{[0,\infty)} e^{-\beta s}\,d\vartheta_s^\delta$$
$$\leq \mathbb{E}\int_{[0,\infty)} e^{-\beta s}h(dY_s^d) + \mathbb{E}\int_{[0,\infty)} e^{-\beta s}h(dY_s^\delta)$$
$$+ h(G\widehat{y}_0)\mathbb{E}\int_{[0,\infty)} e^{-\beta s}\,d\vartheta_s^\delta.$$

By (8.9) and (8.15),

(8.19) $\vartheta_t^\delta \leq c_1 \sup_{0 \leq s \leq t}|GY^\delta(t) - GY^c(t)| \leq 2c_1 \sup_{0 \leq s \leq t}|GY^c(s)| \leq c_4|GY(t)|,$

where (2.2) is used in the last inequality. Combining this observation with (8.8), we obtain $e^{-\beta t}\mathbb{E}\vartheta_t^\delta \to 0$ as $t \to \infty$. Thus, $\mathbb{E}\int_{[0,\infty)} e^{-\beta s}\,d\vartheta_s^\delta = \frac{1}{\beta}\mathbb{E}\int_{[0,\infty)} e^{-\beta s}\vartheta_s^\delta\,ds$. Combining (8.19) and (8.8), we have, by dominated convergence,

(8.20) $$\mathbb{E}\int_{[0,\infty)} e^{-\beta s}\,d\vartheta_s^\delta \to 0 \quad \text{as } \delta \to 0.$$

Next we show that

$$\mathbb{E} \int_{[0,\infty)} e^{-\beta s} h(dY_s^\delta) \leq \mathbb{E} \int_{[0,\infty)} e^{-\beta s} h(dY_s^c). \tag{8.21}$$

Note that

$$\mathbb{E} \int_{[0,\infty)} e^{-\beta s} h(dY_s^\delta) = \mathbb{E} \sum_{j=1}^{\infty} \int_{[j\delta,(j+1)\delta)} e^{-\beta s} h(dY_s^\delta). \tag{8.22}$$

Also, for $j = 1, 2, \ldots$,

$$\int_{[j\delta,(j+1)\delta)} e^{-\beta s} h(dY_s^\delta) = \int_{[j\delta,(j+1)\delta)} e^{-\beta s} h(\dot{Y}_s^\delta) \, ds \leq h(Y_{j,\delta}), \tag{8.23}$$

where $Y_{j,\delta} = e^{-j\beta\delta}(Y^c(j\delta) - Y^c((j-1)\delta))$ and the last inequality follows upon observing that on the interval $[j\delta, (j+1)\delta)$, $\dot{Y}_s^\delta$ equals $\delta^{-1}(Y^c(j\delta) - Y^c((j-1)\delta))$. Using convexity and radial homogeneity of $h$, we see that

$$h(Y^c(j\delta) - Y^c((j-1)\delta)) \leq \int_{[(j-1)\delta,j\delta)} h(Y^{c,\circ}) \, d|Y^c|_s,$$

where $Y^{c,\circ} = dY^c/d|Y^c|$. Thus,

$$\int_{[(j-1)\delta,j\delta)} e^{-\beta s} h(dY_s^c) \geq e^{-\beta j\delta} \int_{[(j-1)\delta,j\delta)} h(dY_s^c)$$
$$\geq e^{-\beta j\delta} h(Y^c(j\delta) - Y^c((j-1)\delta))$$
$$= h(Y_{j,\delta}).$$

Using the above inequality in (8.23) yields that $\int_{[j\delta,(j+1)\delta)} e^{-\beta s} h(dY_s^\delta) \leq \int_{[(j-1)\delta,j\delta)} e^{-\beta s} h(dY_s^c)$. Combining this observation with (8.22) gives (8.21). Using (8.20) and (8.21) in (8.19), we get

$$\limsup_{\delta \to 0} \mathbb{E} \int_{[0,\infty)} e^{-\beta s} h(d\overline{Y}_s^\delta) \leq \mathbb{E} \int_{[0,\infty)} e^{-\beta s} h(dY_s^d)$$
$$+ \mathbb{E} \int_{[0,\infty)} e^{-\beta s} h(dY_s^c) \tag{8.24}$$
$$= \mathbb{E} \int_{[0,\infty)} e^{-\beta s} h(dY_s).$$

Combining (8.25) and (8.17), we obtain (8.14) and the result follows. □

## APPENDIX

LEMMA A.1. *Let $Y_t$ be an $\overline{\mathcal{F}}_t$-adapted RCLL process with increments in $\mathbb{Y}$. One can find an $\overline{\mathcal{F}}_t$-progressively measurable process $Y_t^\circ$ with values in $\mathbb{Y} \cap S^{p-1}$ such that $\int_{[0,t]} Y_s^\circ \, d|Y|_s = Y_t$, $t \geq 0$, a.s.*

PROOF. Since $\overline{\mathcal{F}}_t$ satisfies the usual conditions, it is right continuous. Moreover, since $Y$ is $\overline{\mathcal{F}}_t$-adapted and RCLL, it is $\overline{\mathcal{F}}_t$-progressively measurable (see Proposition 1.13 of [24]). Hence, by Appendix D of [15], there exists a progressively measurable process $Y_t^\circ$ such that $|Y_t^\circ| \leq 1$, $t \geq 0$, and $Y_t = \int_{[0,t]} Y_t^\circ \, d|Y|_t$ a.s. (in [15] $Y_t$ is left continuous and the integral is over $[0,t)$, but the adaptation is clear). It remains to show that

(A.1) $$\mathbb{P}(|Y_t^\circ| = 1, \ d|Y|\text{-a.e.}) = 1.$$

If $y$ is a deterministic RCLL path with increments in $\mathbb{Y}$ and $I$ is a finite interval of $\mathbb{R}_+$,

$$\int_I d|y|_t = \sup \sum_k \left| \int_{I_k} y_t^\circ \, d|y|_t \right| \leq \sup \sum_k \int_{I_k} |y_t^\circ| \, d|y|_t = \int_I |y_t^\circ| \, d|y|_t,$$

where the suprema range over partitions $(I_k)$ of $I$. The above inequality for intervals implies the same for Borel sets $S$ of $\mathbb{R}_+$, thus,

(A.2) $$\int_S d|y|_t \leq \int_S |y_t^\circ| \, d|y|_t.$$

Fix $T > 0$. Given $\varepsilon > 0$, (A.2) implies

$$\int_{[0,T]} \mathbb{1}_{|Y^\circ| < 1-\varepsilon} \, d|Y| \leq \int_{[0,T]} \mathbb{1}_{|Y^\circ| < 1-\varepsilon} |Y^\circ| \, d|Y|$$

$$\leq (1-\varepsilon) \int_{[0,T]} \mathbb{1}_{|Y^\circ| < 1-\varepsilon} \, d|Y|, \qquad \text{a.s.}$$

and, therefore, $E \int_{[0,T]} \mathbb{1}_{|Y^\circ| < 1-\varepsilon} \, d|Y| = 0$. Since $\lim_{\varepsilon \to 0} \int_{[0,T]} \mathbb{1}_{|Y^\circ| < 1-\varepsilon} \, d|Y| = \int_{[0,T]} \mathbb{1}_{|Y^\circ| < 1} d|Y|$ for a.e. $\omega$, we have by Fatou's lemma that $\mathbb{E} \int_{[0,T]} \mathbb{1}_{|Y^\circ| < 1} \, d|Y| = 0$, proving (A.1). $\square$

For the result below, note that $G$ has full row rank, that is, $k$, and $\mathbb{Y}$ has nonempty interior, and therefore, one can find $k$ linearly independent elements of $\mathbb{U}$.

LEMMA A.2. *Let $u_i \in \mathbb{U}$, $i = 1, \ldots, k$, be linearly independent unit vectors. Let $\mathbb{A} = \text{cone}(u_i : i = 1, \ldots, k)$. Then there exists a linear map $\varrho : \mathbb{A} \to \mathbb{Y}$ such that $u = G\varrho(u)$ for all $u \in \mathbb{A}$.*

PROOF. For $i = 1, \ldots, k$, fix $y_i \in \mathbb{Y}$ such that $Gy_i = u_i$. Every $u \in \mathbb{A}$ can be written as $\sum_{i=1}^k \alpha_i u_i$, where $\alpha_i = \alpha_i(u) \geq 0$. Also, $\alpha_i(u)$ are uniquely determined by and depend linearly on $u \in \mathbb{A}$. Setting

$$\varrho(u) = \sum_{i=1}^k \alpha_i(u) \, y_i, \qquad u \in \mathbb{A},$$

yields the result. □

PROOF OF CLAIM (8.3). From Lemma 4.5, we know that $V$ is continuous and so uniformly continuous on $B_{2\varepsilon}(x) \cap \mathbb{X}$. Thus, there exists a $\varsigma \in (0, \infty)$ such that

(A.3)  $|V(x_1) - V(x_2)| < \delta/3$     for all $x_1, x_2 \in B_{2\varepsilon}(x) \cap \mathbb{X}$ s.t. $|x_1 - x_2| \leq \varsigma$.

Let the polyhedral cone $\mathbb{A} \subset \mathbb{U}$ and the linear map $\varrho$ be as in Lemma A.2, and denote by $|\varrho|$ the operator norm of $\varrho$. By making $\varsigma$ smaller if necessary, we can assume that $\varsigma|\varrho| < \delta/3$. Now let $\mathcal{Z}_n \doteq B_{2\varepsilon}(x) \cap \mathbb{X} \cap n^{-1}\mathbb{Z}^k$. Note that $\mathcal{Z}_n$ is a finite set and for $n$ sufficiently large, we can find a measurable map $\vartheta_n : B_\varepsilon(x) \cap \mathbb{X} \to \mathcal{Z}_n$ such that

(A.4)  $\vartheta_n(\xi) - \xi \in \mathbb{A}$   and   $|\vartheta_n(\xi) - \xi| \leq \varsigma$     for all $\xi \in B_\varepsilon(x) \cap \mathbb{X}$.

Fix such $n$ and suppress it from the notation. For each $\xi \in \mathcal{Z}$, find $Y_\xi \in \mathcal{A}(\xi)$ such that

(A.5)    $V(\xi) \geq \mathbb{E}\int_{[0,\infty)} e^{-\beta s}(\ell(X_\xi(s))\,ds + h(dY_\xi(s))) - \delta/3,$

where $X_\xi$ is the controlled process corresponding to $Y_\xi$. Note that for each $\xi \in \mathcal{Z}$, as in the proof of (8.1), we can find measurable function $F_\xi : [0, \infty) \times (\mathbb{R}^k)^{\mathbb{Z}_+} \to \mathbb{Y}$ such that

$$F_\xi(s, \{B(t_i \wedge s)\}_{i \geq 1}) = Y_\xi(s), \qquad s \in [0, \infty).$$

Now define for $s \geq 0$ $Y^*(s) \doteq F_{\vartheta(X(\theta))}(s, \{B((t_i \wedge s) + \theta) - B(\theta)\}_{i \geq 1})$. Finally, define $\widetilde{Y}_s = Y_s \mathbb{1}_{s<\theta} + (Y^*_s + Y_\theta + \varrho(\vartheta(X_\theta) - X_\theta))\mathbb{1}_{s \geq \theta}$ and let $\widetilde{X}$ be the corresponding controlled process. By construction and using (4.6), it is easy to check that

$$\mathbb{E}\left[\int_{[\theta,\infty)} e^{-\beta s}(\ell(\widetilde{X}_s)\,ds + h(d\widetilde{Y}_s)) - e^{-\beta\theta}h(\Delta Y_\theta)\Big|\mathcal{F}_\theta\right]$$
$$\leq e^{-\beta\theta}V(X_\theta) + \delta. \qquad \square$$

**Acknowledgments.** The authors are grateful to Ruth Williams for many valuable suggestions and for pointing out the reference [3]. They would also like to thank Steve Shreve for helpful discussions. Part of the work on this paper was done when the second author was visiting the Technion. Hospitality of the EE Department at the Technion is gratefully acknowledged.

# REFERENCES

[1] BARDI, M. and CAPUZZO-DOLCETTA, I. (1997). *Optimal Control and Viscosity Solutions of Hamilton–Jacobi–Bellman equations*. Birkhäuser, Boston. MR1484411

[2] BENEŠ,V. E., SHEPP, L. A. and WITSENHAUSEN, H. S. (1980/81). Some solvable stochastic control problems. *Stochastics* **4** 39–83. MR0587428

[3] BOHM, V. (1975). On the continuity of the optimal policy set for linear programs. *SIAM J. Appl. Math.* **28** 303–306. MR0371390

[4] BORKAR, V. (1989). *Optimal Control of Diffusion Processes*. Longman Scientific and Technical, Harlow. MR1005532

[5] BRAMSON, M. and WILLIAMS, R. J. (2003). Two workload properties for Brownian networks. *Queueing Systems* **45** 191–221. MR2024178

[6] CAPUZZO-DOLCETTA, I. and LIONS, P. L. (1990). Hamilton–Jacobi equations with state constraints. *Trans. AMS* **318** 643–683. MR0951880

[7] CHOW, P. L., MENALDI, J. L. and ROBIN, M. (1985). Additive control of stochastic linear systems with finite horizon. *SIAM J. Control Optim.* **23** 858–899. MR0809540

[8] CRANDALL, M. G., ISHII, H. and LIONS, P. L. (1992). User's guide to viscosity solutions of second order partial differential equations. *Bull. Amer. Math. Soc.* (*N.S.*) **27** 1–67. MR1118699

[9] CRANDALL, M. G. and LIONS, P. L. (1987). Remarks on the existence and uniqueness of unbounded viscosity solutions of Hamilton–Jacobi equations. *Illinois J. Math.* **31** 665–688. MR0909790

[10] DUFFIE, D., FLEMING, W., SONER, H. M. and ZARIPHOPOULOU, T. (1997). Hedging in incomplete markets with HARA utility. *J. Econom. Dynam. Control* **21** 753–782. MR1455755

[11] DUPUIS, P. and ISHII, H. (1991). On Lipschitz continuity of the solution mapping to the Skorokhod problem, with applications. *Stochastics Stochastics Rep.* **35** 31–62. MR1110990

[12] DUPUIS, P. and RAMANAN, K. (1999). Convex duality and the Skorokhod problem. I, II. *Probab. Theory Related Fields* **2** 153–195, 197–236. MR1720348

[13] DURRETT, R. (1995). *Probability: Theory and Examples*, 2nd ed. Duxbury Press, N. Scituate, MA. MR1609153

[14] EVANS, L. C. (1979). A second-order elliptic equation with gradient constraint. *Comm. Partial Differential Equations* **4** 555–572. MR0529814

[15] FLEMING, W. and SONER, H. (1993). *Controlled Markov Processes and Viscosity Solutions*. Springer, New York. MR1199811

[16] HARRISON, J. M. (1988). Brownian models of queueing networks with heterogeneous customer population. In *Stochastic Differential Systems Stochastic Control Theory and Their Applications* 147–186. Springer, New York. MR0934722

[17] HARRISON, J. M. (2000). Brownian models of open processing networks: Canonical representation of workload. *Ann. Appl. Probab.* **10** 75–103. MR1765204

[18] HARRISON, J. M. and TAKSAR, M. I. (1983). Instantaneous control of Brownian motion. *Math. Oper. Res.* **8** 439–453. MR0716123

[19] HARRISON, J. M. and VAN MIEGHEM, J. A. (1997). Dynamic control of Brownian networks: State space collapse and equivalent workload formulation. *Ann. Appl. Probab.* **7** 747–771. MR1459269

[20] ISHII, H. (1989). On uniqueness and existence of viscosity solutions of fully nonlinear second-order elliptic PDEs. *Comm. Pure Appl. Math.* **42** 15–45. MR0973743

[21] Ishii, H. and Koike, S. (1983). Boundary regularity and uniqueness for an elliptic equation with gradient constraint. *Comm. Partial Differential Equations* **8** 317–346. MR0693645
[22] Ishii, H. and Loreti, P. (2002). A class of optimal control problems with state constraint. *Indiana Univ. Math. J.* **51** 1167–1196. MR1947872
[23] Karatzas, I. (1983). A class of singular stochastic control problems. *Adv. in Appl. Probab.* **15** 225–254. MR0698818
[24] Karatzas, I. and Shreve, S. E. (1991). *Brownian Motion and Stochastic Calculus*, 2nd ed. Springer, New York. MR1121940
[25] Katsoulakis, M. A. (1994). Viscosity solutions of second order fully nonlinear elliptic equations with state constraints. *Indiana Univ. Math. J.* **43** 493–519. MR1291526
[26] Lasry, J. M. and Lions, P. L. (1989). Nonlinear elliptic equations with singular boundary conditions and stochastic control with state constraints. I. The model problem. *Math. Ann.* **283** 583–630. MR0990591
[27] Lions, P. L. (1982). *Generalized solutions of Hamilton–Jacobi Equations*. Pitman, Boston. MR0667669
[28] Martins, L. F., Shreve, S. E. and Soner, H. M. (1996). Heavy traffic convergence of a controlled multi-class queuing system. *SIAM J. Control. Optim.* **34** 2133–2171. MR1416504
[29] Menaldi, J. L. and Robin, M. (1983). On some cheap control problems for diffusion processes. *Trans. AMS* **278** 771–802. MR0701523
[30] Rockafellar, R. T. (1970). *Convex Analysis*. Princeton Univ. Press. MR0274683
[31] Shreve, S. E. and Soner, H. M. (1994). Optimal investment and consumption with transaction costs. *Ann. Appl. Probab.* **4** 609–692. MR1284980
[32] Soner, H. M. (1986). Optimal control with state constraint. I. *SIAM J. Control. Optim.* **24** 552–561. MR0838056
[33] Wong, E. (1971). Representation of martingales, quadratic variation and applications. *SIAM J. Control Optim.* **9** 621–633. MR0307324

Department of Electrical Engineering  
Technion–Israel Institute of Technology  
Haifa 32000  
Israel  
E-mail: atar@ee.technion.ac.il

Department of Statistics  
University of North Carolina  
Chapel Hill, North Carolina 27599-3260  
USA